\documentclass[12pt]{iopart}
\usepackage{iopams}
\usepackage{setstack}
\usepackage{amssymb}
\usepackage{amsfonts}
\usepackage{float}

\usepackage{amsthm}

\usepackage{bbding}
\usepackage[dvips]{graphicx}
\graphicspath{{./figures/}}

\usepackage{subfigure}

\newcommand{\norm}[1]{\left\Vert#1\right\Vert}

\newcommand{\abs}[1]{\left\vert#1\right\vert}
\newcommand{\RE}{{\rm I}\!{\rm R}} 
\newcommand{\Y}{\mathcal{Y}}
\newcommand{\X}{\mathcal{X}}
\newcommand{\D}{\mathcal{D}}
\newcommand{\NU}{\mathcal N}
\newcommand{\R}{\mathcal{R}}
\newcommand{\Z}{\mathcal{Z}}

\newcommand{\ce}{{\scriptscriptstyle 0}}

\newcommand{\s}{\scriptscriptstyle}
\newcommand{\Nu}{\text{\Large$\nu$}}
 \newtheorem{thm}{Theorem}[section]
 \newtheorem{cor}[thm]{Corollary}
 \newtheorem{lem}[thm]{Lemma}
 
 \newtheorem{defn}[thm]{Definition}
 \newtheorem{rem}[thm]{Remark}

\date{\today}
\begin{document}
\title[Solutions of generalized Tikhonov-Phillips functionals]{Existence, uniqueness and stability of solutions of generalized Tikhonov-Phillips functionals}

\author{G. L. Mazzieri $^{1,2}$, R. D. Spies$^{1,3}$ and K. G. Temperini$^{1,4}$}

\address{(1) Instituto de Matem\'{a}tica Aplicada del
Litoral, IMAL, CONICET-UNL, G\"{u}emes 3450, S3000GLN, Santa Fe, Argentina. \\
(2) Departamento de Matem\'{a}tica, Facultad de Bioqu\'{\i}mica y Ciencias Biol\'{o}gicas,
Universidad Nacional del Litoral, Santa Fe, Argentina. \\
(3) Departamento de Matem\'{a}tica, Facultad de Ingenier\'{\i}a Qu\'{\i}mica,
Universidad Nacional del Litoral, Santa Fe, Argentina. \\
(4) Departamento de Matem\'{a}tica, Facultad de Humanidades y
Ciencias, Universidad Nacional del Litoral, Santa Fe, Argentina.}

\eads{\mailto{gmazzieri@hotmail.com}, \mailto{rspies@santafe-conicet.gov.ar},
\mailto{ktemperini@santafe-conicet.gov.ar}}

\begin{abstract}
The Tikhonov-Phillips method is widely used for regularizing
ill-posed inverse problems mainly due to the simplicity of its formulation as an
optimization problem. The use of different penalizers in the
functionals associated to the corresponding optimization problems
has originated a variety other methods which can be considered as
``variants" of the traditional Tikhonov-Phillips method of order
zero. Such is the case for instance of the Tikhonov-Phillips
method of order one, the total variation regularization method,
etc. In this article we find sufficient conditions on the
penalizers in generalized Tikhonov-Phillips functionals which
guarantee existence and uniqueness and stability of the minimizers.
The particular cases in which the penalizers are given by the bounded
variation norm, by powers of seminorms and by linear combinations
of powers of seminorms associated to closed operators, are studied.
Several examples are presented and a few results on image restoration are shown.
\end{abstract}



\maketitle
\section{Introduction}\label{sec:introduction}
In a quite general framework an inverse problem can be formulated
as the need for determining $x$ in an equation of the form
\begin{equation}\label{eq:prob-inv}
Tx=y,
\end{equation}
where $T$ is a linear bounded operator between two infinite
dimensional Hilbert spaces $X$ and $Y$ (in general these will be
function spaces), the range of $T$, $\R(T)$, is non-closed and $y$
is the data, supposed to be known, perhaps with a certain degree
of error. It is well known that under these hypotheses, problem
(\ref{eq:prob-inv}) is ill-posed in the sense of Hadamard
(\cite{ref:Hadamard-1902}). In this case the ill-posedness is a
result of the unboundedness of $T^\dag$, the Moore-Penrose
generalized inverse of $T$. The Moore-Penrose generalized inverse
is a fundamental tool in the treatment of inverse ill-posed
problems and their regularized solutions, mainly due to the fact
that this operator is strongly related to the least-squares
solutions of problem (\ref{eq:prob-inv}). In fact, the
least-squares solution of minimum norm of problem
(\ref{eq:prob-inv}), also known as the best approximate solution,
is $x^\dag\doteq T^\dag y$, which exists if and only if $y\in
\D(T^\dag)=\R(T)\oplus \R(T)^\bot$. Moreover, for $y\in
\D(T^\dag)$, the set of all least-squares solutions of problem
(\ref{eq:prob-inv}) is given by $x^\dag+\NU(T)$, where $\NU(T)$
denotes the null space of the operator $T$.

The unboundedness of $T^\dag$ has as undesired consequence the
fact that small errors or noise in the data $y$ can result in
arbitrarily large errors in the corresponding approximated
solutions (see \cite{ref:Spies-Temperini-2006},
\cite{ref:Seidman-1980}), turning unstable all standard numerical
approximation methods, making them unsuitable for most
applications and inappropriate from any practical point of view.
The so called ``regularization methods" are mathematical tools
designed to restore stability to the inversion process and consist
essentially of parametric families of continuous linear operators
approximating $T^\dag.$ The mathematical theory of regularization
methods is very wide (a comprehensive treatise on the subject can
be found in the book by Engl, Hanke and Neubauer,
\cite{refb:Engl-Hanke-Neubauer-1996}) and it is of great interest
in a broad variety of applications in many areas such as Medicine,
Physics, Geology, Geophysics, Biology, image restaurarion and
processing, etc.

There exist numerous ways of regularizing an ill-posed inverse
problem. Among the most standard and traditional methods we
mention the Tikhonov-Phillips method (\cite{ref:Phillips-1962},
\cite{ref:Tikhonov-1963-SMD-2}, \cite{ref:Tikhonov-1963-SMD-1}),
truncated singular value decomposition (TSVD), Showalter's method,
total variation regularization  (\cite{ref:Acar-Vogel-1994}), etc.
Among all regularization methods, probably the best known and most
commonly and widely used is the Tikhonov-Phillips regularization
method, which was originally proposed by Tikhonov and Phillips in
1962 and 1963 (see \cite{ref:Phillips-1962},
\cite{ref:Tikhonov-1963-SMD-2}, \cite{ref:Tikhonov-1963-SMD-1}).
Although this method can be formalized within a very general
framework by means of spectral theory
(\cite{refb:Engl-Hanke-Neubauer-1996},
\cite{refb:Dautray-Lions-1990}), the widespread of its use is
undoubtedly due to the fact that it can also be formulated in a
very simple way as an optimization problem. In fact, the
regularized solution of problem (\ref{eq:prob-inv}) obtained by
applying Tikhonov-Phillips method is the minimizer $x_\alpha$ of
the functional
\begin{equation}\label{eq:formulacion variacional TP0}
J_\alpha(x)\doteq \norm{Tx-y}^2+\alpha\norm x^2,
\end{equation}
where $\alpha$ is a positive constant known as the regularization
parameter.

The penalizing term $\alpha\norm x^2$ in (\ref{eq:formulacion
variacional TP0}) not only induces stability but it also
determines certain regularity properties of the approximating
regularized solutions $x_\alpha$ and of the corresponding
least-squares solution which they approximate as $\alpha\to 0^+$.
Thus, for instance, it is well known that minimizers of
(\ref{eq:formulacion variacional TP0}) are always ``smooth" and,
for $\alpha\rightarrow 0^+$, they approximate the least-squares
solution of minimum norm of (\ref{eq:prob-inv}), that is
$\lim_{\alpha\to 0^+} x_\alpha=T^\dag y$. This method is more
precisely known as the Tikhonov-Phillips method of order zero.
Choosing other penalizing terms gives rise to different
approximations with different properties, approximating different
least-squares solutions of (\ref{eq:prob-inv}). Thus, for
instance, the use of $\norm{\bigtriangledown x}^2$ as penalizer
instead of $\norm{x}^2$ in (\ref{eq:formulacion variacional TP0})
originates the so called Tikhonov-Phillips method of order one,
the penalizer $\norm{x}_{\scriptscriptstyle \text{BV}}$ (where
$\norm{\cdot}_{\scriptscriptstyle \text{BV}}$ denotes the bounded
variation norm) gives rise to the so called bounded variation
regularization method introduced by Acar and Vogel in 1994
(\cite{ref:Acar-Vogel-1994}), etc. In particular, in the latter
case, the approximating solutions are only forced to be of bounded
variation rather than smooth and they approximate, for
$\alpha\rightarrow 0^+$, the least-squares solution of problem
(\ref{eq:prob-inv}) of minimum $\norm{\cdot}_{\scriptscriptstyle
\text{BV}}$-norm (see \cite{ref:Acar-Vogel-1994}). This method has
been proved to be a good choice, for instance, in certain image
restoration problems in which it is highly desirable to detect and
preserve sharp edges and discontinuities of the original image.

Hence, the penalizing term in (\ref{eq:formulacion variacional
TP0}) is used not only to stabilize the inversion of the ill-posed
problem but also to enforce certain characteristics on the
approximating solutions and on the particular limiting
least-squares solution that they approximate. As a consequence, it
is reasonable to assume that an adequate choice of the penalizing
term, based on {\it a-priori} knowledge about certain
characteristics of the exact solution of problem
(\ref{eq:prob-inv}), will lead to approximated ``regularized"
solutions which will appropriately reflect those characteristics.

With the above considerations in mind, we shall consider
functionals of the form
\begin{equation}\label{eq:prob opt}
J_{\s W,\alpha}(x)\doteq\norm{Tx-y}^2+\alpha
W(x)\hspace{0.5cm}x\in\mathcal D,
\end{equation}
where $W(\cdot)$ is an arbitrary functional with domain $\mathcal
D\subset \X$ and $\alpha$ is a positive constant.

The purpose of this article is to find sufficient conditions on the
penalizers in generalized Tikhonov-Phillips functionals of the form (\ref{eq:prob opt})
which guarantee existence and uniqueness and stability of the minimizers.
The particular cases in which the penalizers are given by the bounded
variation norm, by powers of seminorms and by linear combinations
of powers of seminorms associated to closed operators, are studied.
Several examples are presented and a few results on image restoration are shown.

\section{Existence and uniqueness for general penalizing terms}\label{sec:existence}

In this section we shall consider the problem of finding
conditions on the penalizer $W(\cdot)$ which guarantee existence
and uniqueness of global minimizers of (\ref{eq:prob opt}).
Previously we will need to introduce a few definitions.
\smallskip

\begin{defn}\label{def:W acotado}
Let $\X$ be a vector space, $W$ a functional defined over a set
$\D\subset\X$ and $A$ a subset of $\D$. We say that $A$ is
$W$-bounded if there exists a constant $k<\infty$ such that
$\abs{W(a)}\leq k$ for every $a\in A.$
\end{defn}
\smallskip

\begin{defn}\label{def:W coercitivo}
($W$-coercivity) Let $\X$ be a vector space and $W$, $F$ two
functionals defined on a set $\D\subset\X$. We say that the
functional $F$ is $W$-coercive if
$\displaystyle\lim_{n\rightarrow\infty}F(x_n)=+\infty$ for every
sequence $\{x_n\}\subset \D$ for which $\lim_{n\rightarrow\infty}
W(x_n)=+\infty.$
\end{defn}

\medskip

\begin{rem}\label{obs: w-coerc implica acot} Note that if the
functional $F$ is $W$-coercive and $W$ is bounded from below, then
all lower level sets for $F$, i.e. all sets of the form $\{x\in\D:
F(x)\leq a\}$ with $a\in\RE$ are $W$-bounded sets.
\end{rem}
\begin{defn}\label{def:W-sls-swls}
Let $\X$ be a normed vector space, $W$, $F$ two functionals with
$Dom(F)\subset Dom(W)\subset\X$. We say that $F$ is
$W$-subsequentially (weakly) lower semicontinuous if for every
$W$-bounded sequence $\{x_n\}\subset Dom(F)$ such that
$x_n\overset{(w)}{\rightarrow} x\in Dom (F)$, there exists a
subsequence $\{x_{n_j}\}\subset\{x_n\}$ such that $F(x)\le
\liminf_{j\to\infty}F(x_{n_j})$. If $F$ is $W$-subsequentially
lower semicontinuous we will simply say that $F$ is $W$-sls.
Similarly, if $F$ is $W$-subsequentially weakly lower
semicontinuous we will say that $F$ is $W$-swls.
\end{defn}

In the following theorem, sufficient conditions on the operator
$T$ and on the functional $W$ guaranteeing the existence and
uniqueness of the minimizer of the functional (\ref{eq:prob opt})
are established.

\medskip
\begin{thm}\label{teo:existencia-1}
(Existence and uniqueness) Let $\X,\,\Y$ be normed vector spaces,
$T\in\mathcal L(\X,\Y)$, $y\in\Y,$ $\D\subset \X$ a convex set and
$W:\mathcal{D}\longrightarrow \mathbb R$ a functional bounded from
below, $W$-subsequentially weakly lower semicontinuous, and such
that $W$-bounded sets are relatively weakly compact in $\X$. More
precisely, suppose that $W$ satisfies the following hypotheses:
\begin{itemize}
\item \textit{(H1): }$\exists \,\,\gamma\geq 0$ such that $W(x)\geq
-\gamma\hspace{0.3cm} \forall\,x\in\mathcal D$.
\item \textit{(H2):\;}for every $W$-bounded sequence $\{x_n\}\subset\D$  such that
$x_{n}\overset{w}{\rightarrow} x\in\D,$ there exists a subsequence
$\{x_{n_j}\}\subset \{x_n\}$ such that $
W(x)\leq\liminf_{j\rightarrow\infty}W(x_{n_j})$.
\item \textit{(H3):} for every $W$-bounded sequence
$\{x_n\} \subset \mathcal D$ there exist a subsequence
$\{x_{n_j}\}\subset \{x_n\}$ and $x\in \D$ such that
$x_{n_j}\overset{w}{\rightarrow} x$.
\end{itemize}
Then the functional $J_{W,\alpha}(\cdot)$ in (\ref{eq:prob opt})
has a global minimizer. If moreover $W$ is convex and $T$ is
injective or $W$ is strictly convex, then such a minimizer is
unique.
\end{thm}
\begin{proof}
First we note that for every sequence $\{z_n\}\subset\D$ we have
that
\begin{equation}\label{eq:dsi normas}
z_n\overset{w}{\longrightarrow} z \Longrightarrow \norm{Tz-y}^2\leq\liminf_{n\rightarrow\infty}\norm{Tz_{n}-y}^2.
\end{equation}
This follows immediately from the continuity of $T$ and the weak
lower semicontinuity of the norm.

Let now  $\{x_n\}\subset \D$ be such that
\begin{equation}\label{eq:Jmin finito}
J_{\s W, \alpha} (x_n)\rightarrow \inf_{x\in\D}J_{\s W, \alpha}(x)\doteq J_{\text{min}}.
\end{equation}
Hypothesis \textit{(H1)} guarantees that
$-\infty<J_{\text{min}}<+\infty.$ From the definition of $J_{\s W,
\alpha}(\cdot)$ and since $\alpha>0$ it follows that $J_{\s W,
\alpha}(\cdot)$ is $W$-coercive.
Suppose now that the sequence $\{x_n\}$ is not $W$-bounded. Then,
there exists a subsequence $\{x_{n_j}\}\subset\{x_n\}$ such that
$W(x_{n_j})\to \infty$, from which, by virtue of the
$W$-coercivity of  $J_{\s W, \alpha}(\cdot)$ it follows that
$J_{\s W, \alpha}(x_{n_j})\to \infty$. This contradicts
(\ref{eq:Jmin finito}). Thus the sequence $\{x_n\}$ is
$W$-bounded. It then follows by hypothesis \textit{(H3)} that
there must exist a sequence $\{x_{n_j}\}\subset \{x_n\}$ and
$\bar{x}\in \D$ such that $x_{n_j}\overset{w}{\rightarrow}
\bar{x}$ and since $W$ satisfies \textit{(H2)} there exists a
subsequence $\{x_{n_{j_k}}\}$ of $\{x_{n_j}\}$ such that
\begin{equation}\label{eq:dsi en teo estabilidad}
W(\bar x)\leq \liminf_{k\rightarrow\infty}W(x_{n_{j_k}}).
\end{equation}
Then 
\begin{eqnarray*}
J_{\s W, \alpha}(\bar{x})&= \norm{T\bar{x}-y}^2+\alpha W(\bar{x})\\
&\leq \liminf_{k\rightarrow\infty}\norm{Tx_{n_{j_k}}-y}^2+\alpha
\liminf_{k\rightarrow\infty}W(x_{n_{j_k}})\quad
(\text{by (\ref{eq:dsi normas}) and (\ref{eq:dsi en teo estabilidad})})\\
&\leq
\liminf_{k\rightarrow\infty}\left(\norm{Tx_{n_{j_k}}-y}^2+\alpha
W(x_{n_{j_k}})\right)\qquad (\text{by prop. of
liminf})\\
&=\liminf_{k\rightarrow\infty}J_{\s W, \alpha}(x_{n_{j_k}})\qquad \qquad \qquad \qquad \qquad (\text{by def. of }J_{\s W, \alpha})\\
&=\lim_{n\rightarrow\infty}J_{\s W, \alpha}(x_n) \quad (\text{by
(\ref{eq:Jmin finito}) and since $\{x_{n_{j_k}}\}$ is a
subseq. of  $\{x_n\}$} ) \\
&=J_{\text{min}}.
\end{eqnarray*}
It then follows that $J_{\s W, \alpha}(\bar{x})=J_{\text{min}}$.
This proves the existence of a global minimizer of (\ref{eq:prob
opt}). For the uniqueness, note that under the hypothesis that $W$
be convex and $T$ be injective or $W$ be strictly convex, one has
that the functional $J_{\s W, \alpha}(\cdot)$ is strictly convex
and therefore the global minimizer is unique. \hfill
\end{proof}

\medskip

\begin{rem}\label{obs:cambio de hip en existencia} Note that in the previous theorem
the convexity of\, $\D$ is not needed for the existence. Note also
that if we replace hypotheses \textit{(H2)} and \textit{(H3)} on
the functional $W$ by the assumptions that $W$ be $W$-sls and that
$W$-bounded sets be relatively compact in $\X$, i.e. by the
following hypotheses:
\begin{itemize}
\item \textit{(H2'):} for every $W$-bounded sequence $\{x_n\}\subset\D$ such that
$x_{n}\rightarrow x\in\D,$ there exists a subsequence
$\{x_{n_j}\}\subset \{x_n\}$ such that
$W(x)\leq\liminf_{j\rightarrow\infty}W(x_{n_j});$
\item \textit{(H3'):} for every $W$-bounded sequence $\{x_n\}
\subset \D$ there exist a subsequence $\{x_{n_j}\}\subset \{x_n\}$
and $x\in \D$ such that $x_{n_j}\rightarrow x$,
\end{itemize}
then both existence and uniqueness remain valid.
\end{rem}

\medskip

\begin{rem}\label{remark:xxx} Note that hypothesis \textit{(H3')} is stronger than
\textit{(H3)} which, in turn, is stronger than the hypothesis that
every $W$-bounded set be weakly precompact. Also, \textit{(H2')}
is weaker than \textit{(H2)} which in turn is weaker than the
hypothesis that $W$ be weakly lower semicontinuous.
\end{rem}

\medskip

\begin{rem}\label{remark:W-is-a-norm} If $\X$ is a reflexive Banach
space and $W(\cdot)$ is a norm defined on a subspace $\D$ of $\X$,
which is on $\D$ equivalent or stronger that the norm of $X$, then
it follows that $W$ satisfies hypothesis \textit{(H1)},
\textit{(H2)}, \textit{(H3)} and therefore the functional
(\ref{eq:prob opt}) has a global minimizer on $\D$. If moreover
$T$ is injective or the normed space $\left(\D,W(\cdot)\,\right)$
is complete and separable or Hilbert, then such a minimizer is
unique.
\end{rem}

Observe that hypothesis \textit{(H1)}, \textit{(H2)} and
\textit{(H3)} as well as \textit{(H2')} and \textit{(H3')} impose
conditions only on the penalizer $W(\cdot)$ and not on $T$, so
that the corresponding existence and uniqueness results hold for
any bounded linear operator $T$. It is therefore not surprising
that those conditions can be relaxed if some information on $T$ in
connection to $W(\cdot)$ is provided. The next theorem shows a
result in this direction.
\begin{thm}\label{teo:existence-weak} Let $\X$, $\Y$ be normed
spaces, $T\in\cal{L}(\X,\Y)$, $\D\subset \X$ a convex set and $W$
a real functional on $\D$. Consider the following standing
hypotheses:
\begin{itemize}
\item \textit{(I2):} $W$ is $T$-$W$-swls, i.e for every sequence
$\{x_n\}\subset \D$ such that $\{\|Tx_n\|+W(x_n)\}$ is bounded in
$\mathbb{R}$ (in the sequel we shall refer to such a sequence as a
``$T$-$W$ bounded sequence'') and $x_{n}\overset{w}{\rightarrow}
x\in\D$, there exists a subsequence $\{x_{n_j}\}\subset \{x_n\}$
such that $W(x)\leq\liminf_{j\rightarrow\infty}W(x_{n_j})$.
\item \textit{(I3):} $T$-$W$-bounded sets are relatively weakly compact in $\X$, i.e., for every $T$-$W$-bounded sequence
$\{x_n\} \subset \mathcal D$ there exist a subsequence
$\{x_{n_j}\}\subset \{x_n\}$ and $x\in \D$ such that
$x_{n_j}\overset{w}{\rightarrow} x$.
\end{itemize}
If $T$ and $W(\cdot)$ satisfy the hypotheses \textit{(H1)},
\textit{(I2)} and \textit{(I3)}, then the functional
$J_{W,\alpha}(\cdot)$ in (\ref{eq:prob opt}) has a global
minimizer. If moreover $W$ is convex and $T$ is injective or $W$
is strictly convex, then such a minimizer is unique.
\end{thm}
\begin{proof}
Let $\{x_n\}$ be a minimizing sequence of $J_{\s
W,\alpha}(\cdot)$. From the definition of $J_{W,\alpha}(\cdot)$ it
follows that $\{x_n\}$ is $T$-$W$-bounded. Then by \textit{(I3)}
there must exist $\{x_{n_j}\}\subset\{x_n\}$ and $\bar x\in \D$
such that $x_{n_j}\overset{w}{\rightarrow} \bar x$. Now by virtue
of \textit{(I2)} there exists $\{x_{n_{j_k}}\}\subset\{x_{n_j}\}$
such that $W(\bar x)\le \liminf_{k\to\infty} W(x_{n_{j_k}})$.
Following now the same steps as in Theorem \ref{teo:existencia-1}
we obtain that
$$
J_{\s W, \alpha}(\bar x)=\min_{x\in \D} J_{\s W, \alpha}(x).
$$
If $W$ is convex and $T$ is injective or $W(\cdot)$ is strictly
convex, uniqueness follows from the strict convexity of $J_{\s W,
\alpha}(\cdot)$ on $\D$.
\end{proof}

\begin{rem}\label{obs:H-weakers-1} Note that hypotheses \textit{(I2)} and
\textit{(I3)} are weaker that \textit{(H2)} and \textit{(H3)},
respectively. Also note that both \textit{(I2)} and \textit{(I3)}
hold, for instance if $\X$ is reflexive, $W(\cdot)$ is
subsequentially weakly lower semicontinuous and $T$ and $W$ are
complemented, i.e. there exists a positive constant $c$ such that
$\|Tx\|^2+W(x)\ge c\,\|x\|^2\; \forall\, x \in \D$.
\end{rem}

\begin{rem}\label{obs:H-weakers-2} Just like in Theorem
\ref{teo:existencia-1}, in Theorem \ref{teo:existence-weak}, the
convexity of $\D$ is not needed for the existence. Also note that
if hypothesis \textit{(I2)} and \textit{(I3)} are replaced by the
assumption that $W$ be $T$-$W$-sls and that $T$-$W$-bounded sets
be relatively compact in $\X$, i.e. by the hypoteses:
\begin{itemize}
\item \textit{(I2'):} for every $T$-$W$-bounded sequence $\{x_n\}\subset\D$ such that
$x_{n}\rightarrow x\in\D,$ there exists a subsequence
$\{x_{n_j}\}\subset \{x_n\}$ such that
$W(x)\leq\liminf_{j\rightarrow\infty}W(x_{n_j});$
\item \textit{(I3'):} for every $T$-$W$-bounded sequence $\{x_n\}
\subset \D$ there exist a subsequence $\{x_{n_j}\}\subset \{x_n\}$
and $x\in \D$ such that $x_{n_j}\rightarrow x$,
\end{itemize}
then the results of Theorem \ref{teo:existence-weak} remain valid.
\end{rem}
%


\section{Stability}\label{sec:stability}

As it was previously mentioned, inverse ill-posed problems appear
in a wide variety of applications in diverse areas. Solving these
problems usually involves several steps starting from modeling,
through measurements and data acquisition for the experiment under
study, to the discretization of the mathematical model and the
derivation of numerical approximations for the regularized
solutions. All these steps entail intrinsic errors, many of which
are unavoidable. For this reason, in the context of the study of
inverse ill-posed problems from the optic of Tikhonov-Phillips
methods with general penalizing terms, it is of particular
interest to analyze the stability of the minimizers of the
functional (\ref{eq:prob opt}) under different types of
perturbations. To proceed with some results in this direction we
shall need the following definitions.
\medskip

\begin{defn}\label{def:coercitivity}($W$-coercivity) Let $\X$ be a vector space, $W,$
$F_n,\,n=1,2,\ldots$, functionals defined on a set $\D\subset \X$.
We will say that the sequence $\{F_n\}$ is $W$-coercive if
$\lim_{n\rightarrow\infty}F_n(x_n)=+\infty$ for every sequence
$\{x_n\}\subset\D$ for which $\lim_{n\rightarrow\infty}
W(x_n)=+\infty.$
\end{defn}

\smallskip

\begin{defn} \label{def:consistency} (consistency)
Let $\X$ be a vector space and $W, F, F_n, n=1,2,...,$ functionals
defined on a set $\D\subset \X$. We will say that the sequence
$\{F_n\}$ is consistent for $F$ if $F_nx\rightarrow Fx$ for every
$x\in\D$. We will say that the sequence $\{F_n\}$ is $W$-uniformly
consistent for $F$ if  $F_nx\rightarrow Fx$ uniformly on every
$W$-bounded set, that is if for any given $c>0$ and $\epsilon>0$,
there exists $N=N(c,\epsilon)$ such that $|F_n(x)-F(x)|<\epsilon$
for every $n\geq N$ and every $x\in \D$ such that $\abs{W(x)}\leq
c$.
\end{defn}
In the following theorem we present a weak stability result for
the minimizers of a general functional on a normed space.

\medskip
\begin{thm}\label{teo:estabilidad} Let $\X$ be a normed vector space,
$\D$ a subset of $\X$, $W:\mathcal{D} \longrightarrow \mathbb R$ a
functional satisfying the hypotheses \textit{(H1)} and
\textit{(H3)} of Theorem \ref{teo:existencia-1} (i.e. there exists
$\gamma>0$ such that $W(x)\geq -\gamma$ for every $x\in \D$ and
every $W$-bounded sequence contains a weakly convergent
subsequence with limit in $\D$), $J,J_n,\, n=1,2,\ldots,$
functionals on $\mathcal D$ such that $J$ is $W$-swls and
$\{J_n\}$ is $W$-coercive and $W$-uniformly consistent for $J$.
Assume further that there exists a unique global minimizer $\bar
x\in\D$ of $J$ and that each functional $J_n$ also possesses on
$\D$ a global minimizer $x_n$ (not necessarily unique). Then
$x_n\overset{w}{\rightarrow} \bar x$.
\end{thm}

\begin{proof}
Since for each  $n\in\mathbb N$, $x_n$ minimizes the functional
$J_n$ we have that $J_n(x_n)\leq J_n(\bar x)$. Then
\begin{eqnarray}\label{eq:1 teo estabilidad}
\limsup_{n\rightarrow\infty} J_n(x_n)&\leq
\limsup_{n\rightarrow\infty} J_n(\bar x)= J(\bar x)<\infty,
\end{eqnarray}
where the equality follows from the hypothesis that the sequence
$\{J_n\}$ is $W$-uniformly consistent for $J$. From (\ref{eq:1 teo
estabilidad}), the hypothesis \textit{(H1)} on $W$ and the
hypothesis of $W$-coercitivity of  $\{J_n\}$ it then follows that
the sequence $\{x_n\}$ is $W$-bounded.

Suppose now that the sequence $\{x_n\}$ does not converge weakly
to $\bar x$. Then there exists a subsequence $\big\{x_{n_j}\big\}$
of $\{x_n\}$ such that no subsequence $\big\{x_{n_{j_k}}\big\}$ of
$\left\{x_{n_j}\right\}$ converges weakly to $\bar x$. On the
other hand, since the sequence $\left\{x_{n_j}\right\}$ is
$W$-bounded (since the original sequence is) hypothesis
(\textit{H3}) on $W$ implies that there exist $x^*\in\D$ and a
subsequence $\left\{x_{n_{j_k}}\right\}$ of
$\left\{x_{n_{j}}\right\}$ such that
$x_{n_{j_k}}\overset{w}{\rightarrow}x^*$. It then follows that
$x^*\neq\bar x$.

On the other hand, since the sequence $\{x_{n_{j_k}}\}$ is
$W$-bounded, $x_{n_{j_k}}\overset{w}{\rightarrow}x^*$ and $J$ is
$W$-swls, it follows that there exists a subsequence
$\{x_{n_{j_{k_\ell}}}\}\subset \{x_{n_{j_k}}\}$ such that
\begin{eqnarray}\label{eq:numerar1}
J(x^*)&\leq \liminf_{\ell\rightarrow\infty}
J\Big(x_{n_{j_{k_\ell}}}\Big).
\end{eqnarray}
Also, since the sequence $\big\{x_{n_{j_{k_\ell}}}\big\}$ is
$W$-bounded and $\big\{J_n\big\}$ is $W$-uniformly consistent for
$J$, it follows that
\begin{eqnarray}\label{eq:numerar2}
\lim_{\ell\rightarrow\infty}\left(J\left(x_{n_{j_{k_\ell}}}\right)-J_{n_{j_{k_\ell}}}\left(x_{n_{j_{k_\ell}}}\right)\right)=0.
\end{eqnarray}
Hence
\begin{eqnarray*} 
J(x^*)&\leq \liminf_{\ell\rightarrow\infty} J\left(x_{n_{j_{k_\ell}}}\right) & \hspace{-2cm}(\text{by (\ref{eq:numerar1})})\\
&\leq \limsup_{\ell\rightarrow\infty} J\left(x_{n_{j_{k_\ell}}}\right) & \\
&=  \limsup_{\ell\rightarrow\infty}\left[\left( J\left(x_{n_{j_{k_\ell}}}\right)-J_{n_{j_{k_\ell}}}\left(x_{n_{j_{k_\ell}}}\right) \right) +J_{n_{j_{k_\ell}}}\left(x_{n_{j_{k_\ell}}}\right)\right]\\
&\leq
\limsup_{\ell\rightarrow\infty}\left(J\left(x_{n_{j_{k_\ell}}}\right)-J_{n_{j_{k_\ell}}}\left(x_{n_{j_{k_\ell}}}\right)\right)+
\limsup_{\ell\rightarrow\infty} J_{n_{j_{k_\ell}}}\left(x_{n_{j_{k_\ell}}}\right)&\\
&=
\lim_{\ell\rightarrow\infty}\left(J\left(x_{n_{j_{k_\ell}}}\right)-J_{n_{j_{k_\ell}}}\left(x_{n_{j_{k_\ell}}}\right)\right)+
\limsup_{\ell\rightarrow\infty} J_{n_{j_{k_\ell}}}\left(x_{n_{j_{k_\ell}}}\right)&  \hspace{-.3cm}(\text{by (\ref{eq:numerar2})})\\
&=\limsup_{\ell\rightarrow\infty} J_{n_{j_{k_\ell}}}\left(x_{n_{j_{k_\ell}}}\right)& \hspace{-4.5cm}(\text{by (\ref{eq:numerar2})}) \\
&\leq J(\bar x). &\hspace{-4.5cm}(\text{by (\ref{eq:1 teo
estabilidad}), since $\{x_{n_{j_{k_\ell}}}\}\subset\{x_n\}$)}
\end{eqnarray*}
Since $\bar x$ is the unique minimizer of $J$ it follows that
$x^*=\bar x$, contradicting our previous result that $x^*\ne\bar
x$. This contradiction came from assuming that the sequence
$\{x_n\}$ did not converge weakly to $\bar x$. Hence
$x_n\overset{w}{\rightarrow} \bar x$ as we wanted to show.
\end{proof}

Note that by virtue of Remark \ref{remark:xxx}, the hypothesis
that $J$ be $W$-swls in the previous theorem can be replaced by
the hypothesis that $J$ be weakly lower semicontinuous on
$\mathcal D$.

\medskip
In the particular case in which the functionals $J$ and $J_n$ are
of Tikhonov-Phillips type, under certain general conditions on the
penalizer $W(\cdot)$, the previous theorem yields a weak stability
result for the minimizers of the functional (\ref{eq:prob opt}).
In fact we have the following corollary.

\begin{cor} \label{corollary:weak}
Let $\mathcal X$ be a normed vector space, $\mathcal Y$ an inner
product space, $T,\;T_n\in\mathcal L(\mathcal X, \mathcal Y)$,
$n=1,2,\cdots$, $y\in \mathcal Y$, $\alpha>0$, $\mathcal D$ a
subset of $\mathcal X$, $W:\mathcal D\rightarrow \mathbb R$ a
functional satisfying hypotheses (H1), (H2) and (H3) of Theorem
\ref{teo:existencia-1}, $J, J_n, n=1,2,...$, functionals on
$\mathcal D$ defined as follows:
\begin{eqnarray}
J(x)&\doteq\|Tx-y\|^2+\alpha W(x),\label{eq:fun J}\\
J_n(x)&\doteq\|T_nx-y_n\|^2+\alpha_n W(x),\label{eq:fun Jn}
\end{eqnarray}
such that as $n\rightarrow\infty$,\,\,$\alpha_n\rightarrow\alpha$,
$y_n\rightarrow y$ and $T_n x\rightarrow Tx$ uniformly for $x$ in
$W$-bounded sets (i.e. $\{T_n\}$ is $W$-uniformly consistent for
$T$). Suppose further that $J$ has a unique global minimizer $\bar
x$. If $x_n$ is a global minimizer of $J_n$ then
$x_n\overset{w}{\rightarrow} \bar x$.
\end{cor}
\begin{proof}
To prove this corollary it suffices to verify that the functionals
$J$ and $J_n$ satisfy the hypotheses of Theorem
\ref{teo:estabilidad}, that is, verify that $J$ is $W$-swls and
that the sequence $\{J_n\}$ is $W$-coercive and $W$-uniformly
consistent for $J$.

To prove that $J$ is $W$-swls, let $\{x_n\}\subset \mathcal D$ be
a $W$-bounded sequence such that
$x_n\overset{w}{\longrightarrow}x\in\mathcal D$. From the
continuity of $T$ and the weak lower semicontinuity of every norm,
it follow immediately that
\begin{eqnarray}\label{eq:J dsi}
\|Tx-y\|^2\leq \liminf_{n\to\infty} \|Tx_n-y\|^2.
\end{eqnarray}
On the other hand, by \textit{(H2)} it follows that there exists a
subsequence $\{x_{n_j}\}\subset\{x_n\}$ such that
\begin{eqnarray}\label{eq:W dsi}
W(x)\leq \liminf_{j\to\infty} W(x_{n_j}).
\end{eqnarray}
Then,
\begin{eqnarray*}
J(x)&=\|Tx-y\|^2+\alpha W(x)&\\
&\leq\liminf_{j\to\infty} \|Tx_{n_j}-y\|^2+\liminf_{j\to\infty}
\alpha W(x_{n_j})& \text{(by
(\ref{eq:J dsi}) y (\ref{eq:W dsi}))}\\
&\leq\liminf_{j\to\infty}\{\|Tx_{n_j}-y\|^2+\alpha W(x_{n_j})\}&\text{(by property of $\liminf$)}\\
&=\liminf_{j\to\infty} J(x_{n_j}).
\end{eqnarray*}
Hence $J$ is $W$-swls.

Now we will prove that the sequence $\{J_n\}$ is $W$-coercive. For
that let $\{x_n\}\subset \mathcal D$ such that $W(x_n)\rightarrow
+\infty$. Observe that
\begin{eqnarray}\label{eq:Jn coerc}
J_n(x_n)= \|T_nx_n-y_n\|^2+\alpha_n W(x_n)\geq \alpha_n W(x_n).
\end{eqnarray}
Since $W$ satisfies \textit{(H1)} and $\alpha_n\to\alpha>0$, it
follows immediately from (\ref{eq:Jn coerc}) that
$J_n(x_n)\rightarrow +\infty$. Hence $\{J_n\}$ is uniformly
$W$-coercive.

Finally we will show that $\{J_n\}$ is $W$-uniformly consistent
for $J$. For that let $M\subset \D$ be a $W$-bounded set. Since
$\{T_n\}$ is $W$-uniformly consistent for $T$ we have that
$T_nx\rightarrow Tx$ uniformly on $M$ and since $y_n\rightarrow
y$, it follows that $\|T_nx-y_n\|^2\rightarrow \|Tx-y\|^2$
uniformly on $M$. Finally, since
\begin{eqnarray}
|J_n(x)-J(x)|&=\big|\,\|T_nx-y_n\|^2+\alpha_n
W(x)-\|Tx-y\|^2-\alpha
W(x)\big|&\nonumber\\
&\leq
\big|\,\|T_nx-y_n\|^2-\|Tx-y\|^2\big|+\big|(\alpha_n-\alpha)\big|\,
\big|W(x)\big|, \label{eq:no se}&
\end{eqnarray}
it follows that $J_n(x)\rightarrow J(x)$ uniformly for $x\in M$.
Thus $\{J_n\}$ is $W$-uniformly consistent for $J$.

Since $J$ and $\{J_n\}$ satisfy the hypotheses of Theorem
\ref{teo:estabilidad}, the corollary then follows.
\end{proof}

\medskip

\begin{rem}\label{remark:W-is-a-norm2} Note that by virtue of
Remark \ref{remark:W-is-a-norm}, the weak stability result of
Corollary \ref{corollary:weak} holds if \textit{i)} $\X$ is a
reflexive Banach space, \textit{ii)} the penalizer $W(\cdot)$ in
(\ref{eq:fun J}) is a norm defined on a subspace $\D$ of $\X$
which is on $\D$ equivalent or stronger than the original norm in
$\X$ and \textit{iii)} $T$ is injective or the space
$\left(\D,W(\cdot)\,\right)$ is a separable Banach space or a
Hilbert space.
\end{rem}

Hypotheses on Theorem \ref{teo:estabilidad} and Corollary
\ref{corollary:weak} can be weakened if adequate information on
the operator $T$ is available. Before we proceed to the statements
of the corresponding results, we shall need the following
definitions.
\medskip

\begin{defn}\label{def:T-W-coercitivity}($T$-$W$-coercivity) Let $\X$, $\Y$ be vector spaces,
$T\in\cal{L}(\X,\Y)$, $W,$ $F_n,\,n=1,2,\ldots,$ functionals
defined on a set $\D\subset \X$. We will say that the sequence
$\{F_n\}$ is $T$-$W$-coercive if\;
$\lim_{n\rightarrow\infty}F_n(x_n)=+\infty$ for every sequence
$\{x_n\}\subset\D$ for which $\lim_{n\rightarrow\infty}
\|Tx_n\|+W(x_n)=+\infty$.
\end{defn}

\smallskip

\begin{defn} \label{def:T-W-consistency} ($T$-$W$-uniform consistency)
Let $\X$, $\Y$ be vector spaces, $T\in\cal{L}(\X,\Y)$ and $W, F,
F_n, n=1,2,...,$ functionals defined on a set $\D\subset \X$. We
will say that the sequence $\{F_n\}$ is $T$-$W$-uniformly
consistent for $F$ if $F_n\rightarrow F$ uniformly on every
$T$-$W$-bounded set, that is if for any given $c>0$ and
$\epsilon>0$, there exists $N=N(c,\epsilon)$ such that
$|F_n(x)-F(x)|<\epsilon$ for every $n\geq N$ and every $x\in \D$
such that $\|Tx\|+\abs{W(x)}\leq c$.
\end{defn}
\medskip
\begin{thm}\label{teo:estabilidad-weak-2} Let $\X$, $\Y$ be normed vector spaces,
$T\in\cal{L}(\X,\Y)$, $\D$ a subset of $\X$, $W:\mathcal{D}
\longrightarrow \mathbb R$ a functional satisfying hypotheses
\textit{(H1)} of Theorem \ref{teo:existencia-1} and \textit{(I3)}
of Theorem \ref{teo:existence-weak} (i.e. there exists $\gamma>0$
such that $W(x)\geq -\gamma$ for every $x\in \D$ and every
$T$-$W$-bounded sequence contains a weakly convergent subsequence
with limit in $\D$), $J,J_n,\, n=1,2,\ldots,$ functionals on
$\mathcal D$ such that $J$ is $T$-$W$-swls and $\{J_n\}$ is
$T$-$W$-coercive and $T$-$W$-uniformly consistent for $J$. Assume
further that there exists a unique global minimizer $\bar x\in\D$
of $J$ and that each functional $J_n$ also possesses on $\D$ a
global minimizer $x_n$ (not necessarily unique). Then
$x_n\overset{w}{\rightarrow} \bar x$.
\end{thm}
\begin{proof} The proof follows like in Theorem
\ref{teo:estabilidad} with the obvious modifications.
\end{proof}
\begin{cor} \label{corollary:weak-2}
Let $\X$ be normed vector space, $\mathcal Y$ an inner product
space, $T,\;T_n\in\mathcal L(\mathcal X, \mathcal Y)$,
$n=1,2,\cdots$, $y\in \mathcal Y$, $\alpha>0$, $\mathcal D$ a
subset of $\mathcal X$, $W:\mathcal D\rightarrow \mathbb R$ a
functional satisfying hypotheses (H1) of Theorem
\ref{teo:existencia-1} and (I2) and (I3) of Theorem
\ref{teo:existence-weak}, $J, J_n, n=1,2,...$, functionals on
$\mathcal D$ defined as follows:
\begin{eqnarray*}
J(x)&\doteq\|Tx-y\|^2+\alpha W(x),\\ 
J_n(x)&\doteq\|T_nx-y_n\|^2+\alpha_n W(x),
\end{eqnarray*}
such that as $n\rightarrow\infty$,\,\,$\alpha_n\rightarrow\alpha$,
$y_n\rightarrow y$ and $T_n x\rightarrow Tx$ uniformly for $x$ in
$W$-bounded sets (i.e. $\{T_n\}$ is $W$-uniformly consistent for
$T$). Suppose further that $J$ has a unique global minimizer $\bar
x$. If $x_n$ is a global minimizer of $J_n$ then
$x_n\overset{w}{\rightarrow} \bar x$.
\end{cor}
\begin{proof} We will show that $J$ and $\{J_n\}$ satisfy the
hypotheses of Theorem \ref{teo:estabilidad-weak-2}. For that it
suffices to show that $J$ is $T$-$W$-swls and that $\{J_n\}$ is
$T$-$W$ coercive and $T$-$W$-uniformly consistent for $J$. That
$J$ is $T$-$W$-swls follows immediately from \textit{(I2)} and the
weak lower semicontinuity of every norm. The $T$-$W$-uniform
consistency of $\{J_n\}$ for $J$ follows exactly as in the proof
of Corollary \ref{corollary:weak} by noting that $|J_n(x)-J(x)|
\leq\big|\,\|T_nx-y_n\|^2-\|Tx-y\|^2\big|+\big|(\alpha_n-\alpha)\big|\,
\big|W(x)\big|$ and using the fact that $T$-$W$-bounded sets are
also $W$-bounded. Finally, the $T$-$W$-coercivity of $\{J_n\}$
follows easily from the $W$-uniform consistency of $\{T_n\}$ for
$T$.
\end{proof}


\medskip
Next we present a strong stability result for the minimizers of
general functionals on a normed space.

\smallskip

\begin{thm}\label{teo:estabilidad fuerte}
Let $\X$ be a normed vector space, $\D$ a subset of $\X$,
$W:\mathcal{D} \longrightarrow \mathbb R$ a functional satisfying
hypotheses \textit{(H1)} of Theorem \ref{teo:existencia-1} and
\textit{(H3')} of Remark \ref{obs:cambio de hip en existencia}
(i.e., there exists $\gamma>0$ such that $W(x)\geq -\gamma$ for
every $x\in \D$ and every $W$-bounded sequence contains a
convergent subsequence with limit in $\D$), $J,J_n,\,
n=1,2,\ldots,$ functionals on $\mathcal D$ such that $J$ is
$W$-subsequentially lower semicontinuous ($W$-sls) and $\{J_n\}$
is $W$-coercive and $W$-uniformly consistent for $J$. Suppose
further that $J$ has a unique global minimizer $\bar x\in\D$ and
that each functional $J_n$ also possesses on $\D$ a global
minimizer $x_n$ (not necessarily unique). Then $x_n\rightarrow\bar
x$.
\end{thm}

\begin{proof}
For each $n\in\mathbb{N}$, let $x_n$ be a global minimizer of
$J_n$. Following the same steps as those in the proof of Theorem
\ref{teo:estabilidad} it follows that the sequence $\{x_n\}$ is
$W$-bounded.

Suppose now that $\{x_n\}$ does not converge to $\bar x$. Then
there exists a subsequence $\big\{x_{n_j}\big\}$ of $\{x_n\}$ such
that no subsequence of $\left\{x_{n_j}\right\}$ converges to $\bar
x$. On the other hand, since the sequence $\left\{x_{n_j}\right\}$
is $W$-bounded (since the original sequence is), hypothesis
(\textit{H3'}) on the functional $W$ implies that there exist a
subsequence $\left\{x_{n_{j_k}}\right\}$ of
$\left\{x_{n_{j}}\right\}$ and $x^*\in\D$ such that
$x_{n_{j_k}}\rightarrow x^*$. From this it follows that
$x^*\neq\bar x$ and since $J$ is $W$-sls, there exists a
subsequence
$\left\{x_{n_{j_{k_\ell}}}\right\}\subset\left\{x_{n_{j_k}}\right\}$
such that
\begin{eqnarray}\label{eq:numerar1 fuerte}
J(x^*)&\leq \liminf_{\ell\rightarrow\infty}
J\Big(x_{n_{j_{k_\ell}}}\Big).
\end{eqnarray}
Then 
\begin{eqnarray*}
J(x^*)&\leq \liminf_{\ell\rightarrow\infty}
J\left(x_{n_{j_{k_\ell}}}\right)\qquad \qquad\qquad
(\text{by (\ref{eq:numerar1 fuerte})})\\
&\leq\limsup_{\ell\rightarrow\infty}\left[\left( J\left(x_{n_{j_{k_\ell}}}\right)-J_{n_{j_{k_\ell}}}
\left(x_{n_{j_{k_\ell}}}\right) \right)+J_{n_{j_{k_\ell}}}\left(x_{n_{j_{k_\ell}}}\right)\right]\\
&\leq\limsup_{\ell\rightarrow\infty}\left(J\left(x_{n_{j_{k_\ell}}}\right)-J_{n_{j_{k_\ell}}}
\left(x_{n_{j_{k_\ell}}}\right)\right)+\limsup_{\ell\rightarrow\infty}
J_{n_{j_{k_\ell}}}\left(x_{n_{j_{k_\ell}}}\right)\\
&=\limsup_{\ell\rightarrow\infty}
J_{n_{j_{k_\ell}}}\left(x_{n_{j_{k_\ell}}}\right)\quad\text{(since
}\{J_n\}
\text{ is } W\text{-unif. consistent for } J) \\
&\leq \limsup_{n\rightarrow\infty}
J_{n}\left(x_{n}\right)\qquad\qquad (\text{since }\{
x_{n_{j_{k_\ell}}}\}
\subset\{x_n\}) \\
&\leq \limsup_{n\rightarrow\infty} J_{n}\left(\bar x\right)\qquad\qquad(\text{since }x_{n} \text{ minimizes } J_n) \\
&= J(\bar x).\qquad\qquad\qquad(\text{since }\{J_n\} \text{ is }
W\text{-unif. consistent for } J)
\end{eqnarray*}
Hence  $J(x^*)\leq J(\bar x)$ which contradicts the fact that
$\bar x\ne x^*$ and $x^*$ is the unique minimizer of $J$. This
contradiction came from assuming that the sequence $\{x_n\}$ does
not converge to $\bar x$. Hence  $x_n\rightarrow \bar x$. \hfill
\end{proof}
The previous theorem yields a strong stability result for
minimizers of the functional (\ref{eq:prob opt}) in the particular
case in which $J$ and $J_n$ are of Tikhonov-Phillips type. More
precisely we have the following corollary.
\begin{cor} \label{cor:strong-stability-Tikhonov}
Let $\mathcal X$ be a normed vector space, $\mathcal Y$ an inner
product space, $T,\;T_n\in\mathcal L(\mathcal X, \mathcal Y)$,
$n=1,2,\cdots$, $y\in \mathcal Y$, $\alpha>0$, $\mathcal D$ a
subset of $\mathcal X$, $W:\mathcal D\rightarrow \mathbb R$ a
functional satisfying hypotheses (H1) of Theorem
\ref{teo:existencia-1} and (H2') and (H3') of Remark
\ref{obs:cambio de hip en existencia}, $J, J_n, n=1,2,...$,
functionals on $\mathcal D$ defined as follows:
\begin{eqnarray}
J(x)&\doteq\|Tx-y\|^2+\alpha W(x), \\      
J_n(x)&\doteq\|T_nx-y_n\|^2+\alpha_n W(x), 
\end{eqnarray}
such that  as $n\rightarrow\infty$, $\alpha_n\rightarrow\alpha$,
$y_n\rightarrow y$ and $T_n x\rightarrow Tx$ uniformly on
$W$-bounded sets (i.e. $\{T_n\}$ is $W$-uniformly consistent for
$T$). Suppose further that $J$ has a unique global minimizer $\bar
x$. If $x_n$ is a global minimizer of $J_n$ then $x_n\rightarrow
\bar x$.
\end{cor}
\begin{proof}
Since the proof is immediately obtained from Theorem
\ref{teo:estabilidad fuerte} following the same steps as in
Corollary \ref{corollary:weak}, we do not give details here.
\end{proof}
\medskip

Here again, the strong stability results of Theorem
\ref{teo:estabilidad fuerte} and Corollary
\ref{cor:strong-stability-Tikhonov} remain valid under weaker
hypotheses involving both the model operator $T$ and the penalizer
$W$.

\begin{thm}\label{teo:estabilidad fuerte-2}
Let $\X$ be a normed vector space, $\D$ a subset of $\X$,
$W:\mathcal{D} \longrightarrow \mathbb R$ a functional satisfying
hypotheses \textit{(H1)} of Theorem \ref{teo:existencia-1} and
\textit{(I3')} of Remark \ref{obs:H-weakers-2} (i.e., there exists
$\gamma>0$ such that $W(x)\geq -\gamma$ for every $x\in \D$ and
every $T$-$W$-bounded sequence contains a convergent subsequence
with limit in $\D$), $J,J_n,\, n=1,2,\ldots,$ functionals on
$\mathcal D$ such that $J$ is $T$-$W$-subsequentially lower
semicontinuous ($T$-$W$-sls) and $\{J_n\}$ is $T$-$W$-coercive and
$T$-$W$-uniformly consistent for $J$. Suppose further that $J$ has
a unique global minimizer $\bar x\in\D$ and that each functional
$J_n$ also possesses on $\D$ a global minimizer $x_n$ (not
necessarily unique). Then $x_n\rightarrow\bar x$.
\end{thm}
\begin{proof}
The proof of this theorem proceeds exactly as the one of Theorem
\ref{teo:estabilidad fuerte}, by changing the $W$-boundedness,
$W$-sls, $W$-uniform consistency and \textit{(H3')} hypoteses by
$T$-$W$-boundedness, $T$-$W$-sls, $T$-$W$-uniform consistency and
\textit{(I3')}, respectively.
\end{proof}
Here again, the previous strong stability theorem yields a
corresponding stability result for minimizers of the functional
(\ref{eq:prob opt}) in the particular case in which $J$ and $J_n$
are of Tikhonov-Phillips type. This result is given in the
following corollary.
\begin{cor} \label{cor:strong-stability-Tikhonov-2}
Let $\mathcal X$ be a normed vector space, $\mathcal Y$ an inner
product space, $T,\;T_n\in\mathcal L(\mathcal X, \mathcal Y)$,
$n=1,2,\cdots$, $y\in \mathcal Y$, $\alpha>0$, $\mathcal D$ a
subset of $\mathcal X$, $W:\mathcal D\rightarrow \mathbb R$ a
functional satisfying hypotheses (H1) of Theorem
\ref{teo:existencia-1} and (I2') and (I3') of Remark
\ref{obs:H-weakers-2}, $J, J_n, n=1,2,...$, functionals on
$\mathcal D$ defined as follows:
\begin{eqnarray}
J(x)&\doteq\|Tx-y\|^2+\alpha W(x), \\      
J_n(x)&\doteq\|T_nx-y_n\|^2+\alpha_n W(x), 
\end{eqnarray}
such that as $n\rightarrow\infty$, $\alpha_n\rightarrow\alpha$,
$y_n\rightarrow y$ and $T_n x\rightarrow Tx$ uniformly on
$W$-bounded sets (i.e. $\{T_n\}$ is $W$-uniformly consistent for
$T$). Suppose further that $J$ has a unique global minimizer $\bar
x$. If $x_n$ is a global minimizer of $J_n$ then $x_n\rightarrow
\bar x$.
\end{cor}
\begin{proof}
We will show that $J$ and $\{J_n\}$ satisfy the hypotheses of
Theorem \ref{teo:estabilidad fuerte-2}. For that it suffices to
show that $J$ is $T$-$W$-sls and that $\{J_n\}$ is $T$-$W$
coercive and $T$-$W$-uniformly consistent for $J$. The fact that
$J$ is $T$-$W$-sls follows immediately from \textit{(I2')}, the
boundedness of $T$ and the continuity of the norm in $\X$. The
$T$-$W$-uniform consistency of $\{J_n\}$ for $J$ follows exactly
as in the proof of Corollary \ref{corollary:weak} by noting that
$|J_n(x)-J(x)|
\leq\big|\,\|T_nx-y_n\|^2-\|Tx-y\|^2\big|+\big|(\alpha_n-\alpha)\big|\,
\big|W(x)\big|$ and using the fact that $T$-$W$-bounded sets are
also $W$-bounded and the hypothesis of the $W$-uniform consistency
of $\{T_n\}$ for $T$. Finally, also the $T$-$W$-coercivity of
$\{J_n\}$ follows easily from the $W$-uniform consistency of
$\{T_n\}$ for $T$.
\end{proof}

\section{Particular cases} \label{sec:particular-cases}

In this section we present several examples of penalizers
$W(\cdot)$ for which some of the results obtained in the previous
section are valid and therefore, existence, uniqueness and/or
stability for the minimizers of the corresponding generalized
Tikhonov-Phillips functional $J_{\s W,\alpha}(\cdot)$ in
(\ref{eq:prob opt}) are obtained.

\medskip

\subsection {Total variation penalization}\label{ssec:BV}

Bounded variation penalty methods have been studied by Rudin, Osher
and Fatemi in 1992 (\cite{ref:Rudin-Osher-Fatemi-1992}) and Acar and
Vogel in 1994 (\cite{ref:Acar-Vogel-1994}), among others. These
methods have been proved highly successful in certain image
denoising problems where edge preserving is an important issue
(\cite{ref:Chambolle-Lions-1997},
\cite{ref:Chan-Marquina-Mullet-2000}, \cite{ref:Chan-Shen-2002},
\cite{ref:Dobson-Santosa-1996}). Let $d\geq2$, $\Omega\subset\mathbb
R^d$ a convex, bounded set with Lipschitz continuous boundary,
$1\leq p\leq\frac{d}{d-1},\, \X\doteq L^p(\Omega),\, \D\doteq
BV(\Omega),$ where $BV(\Omega)$ denotes the space of functions of
bounded variations on $\Omega$. Recall that $BV(\Omega)=\left\{u\in
L^1(\Omega)\;:\;J_{\s 0}(u)<\infty\right\}$, where $J_{\s
0}(u)\doteq\sup_{v\in\Nu}\int_\Omega(-u\,div\,v)\;dx$ and
$\Nu\doteq\left\{v\in C_{\s 0}^{\s
1}(\Omega;\mathbb{R}^d)\;:\;|v(x)|\le 1 \;\forall\; x\in\Omega
\right\}$ (for $u\in C^1(\Omega)$ one has that $J_{\s
0}(u)=\int_\Omega|\nabla u|\;dx$) and for $u\in BV(\Omega)$ the BV
norm of $u$ is defined by $\norm{u}_{\s BV(\Omega)} \doteq
\norm{u}_{\s L^1(\Omega)}+J_{\s 0}(u)$. Let $W$ be the functional
defined on $\D$ by $W(u)\doteq\norm{u}_{\s BV(\Omega)}$.

We will show that $W(\cdot)$ satisfies the hypotheses
\textit{(H1)}, \textit{(H2)} and \textit{(H3)} of Theorem
\ref{teo:existencia-1}. Clearly $W(\cdot)$ satisfies hypothesis
\textit{(H1)} with $\gamma=0$. Hypothesis \textit{(H3)} follows
immediately from the compact imbedding of $BV(\Omega)$ into
$L^p(\Omega)$ for $1\leq p < \frac{d}{d-1}$ and from the weak
compact imbedding for $p=\frac{d}{d-1}$. These results are
extensions of the Rellich-Kondrachov Theorem and can be found for
example in \cite{refb:Adams-1975} and
\cite{refb:Attouch-Buttazzo-Michaille-2006}. It only remains prove
that $W(\cdot)$ satisfies hypothesis \textit{(H2)}. For that, let
$\{u_n\}\subset \D$ be a $W$-bounded sequence such that
$u_n\overset{w-L^p}{\longrightarrow} u\in\D$. Then,
$u_n\overset{w-L^1}{\longrightarrow} u$ (since $p\geq 1$). From
the weak lower semicontinuity of the $\norm{\cdot}_{L^1(\Omega)}$ norm
and of the functional $J_\ce(\cdot)$ in $L^1(\Omega)$ (see
\cite{ref:Acar-Vogel-1994}), it follows that
\begin{eqnarray}\label{eq:1 ejemplo BV condiciones}
\norm{u}_{L^1(\Omega)}&\leq\liminf_{n\rightarrow\infty}\norm{u_n}_{L^1(\Omega)}\hspace{0.2cm}
\text{and} \hspace{0.2cm}
J_\ce(u)\leq\liminf_{n\rightarrow\infty}J_\ce(u_n).
\end{eqnarray}
Then,
\begin{eqnarray*}
W(u)=\norm{u}_{\s BV(\Omega)}&= \norm{u}_{L^1(\Omega)}+J_\ce(u)&\\
&\leq \liminf_{n\rightarrow\infty}\norm{u_n}_{L^1(\Omega)}+ \liminf_{n\rightarrow\infty}J_\ce(u_n)&\qquad\qquad (\text{by (\ref{eq:1 ejemplo BV condiciones})})\\
&\leq \liminf_{n\rightarrow\infty}\left(\norm{u_n}_{L^1(\Omega)}+ J_\ce(u_n)\right) &\\
&= \liminf_{n\rightarrow\infty}\norm{u_{n}}_{\s BV(\Omega)}&\\
&= \liminf_{n\rightarrow\infty} W(u_{n}),&
\end{eqnarray*}
which proves \textit{(H2)}. Hence $W(\cdot)$ satisfies the
hypotheses of Theorem \ref{teo:existencia-1} and therefore for any
$\alpha>0$, $T\in\mathcal L(\X,\Y)$, ($\Y$ a normed space) the
functional
\begin{equation}\label{BV-functional}
J_{{\s \|\cdot\|_{\s BV}},\,\alpha}(u)\doteq\|Tu-v\|^2 +
\alpha\|u\|_{\s BV(\Omega)}
\end{equation}
has a global minimizer on $BV(\Omega)$. If $T$ is injective then
such a global minimizer is unique. If $T$ is not injective
uniqueness cannot be guaranteed since the $\|\cdot\|_{\s BV}$-norm
is not strictly convex. Also, if $p<\frac{d}{d-1}$ and $J_{{\s
\|\cdot\|_{\s BV}},\,\alpha}(\cdot)$ has a unique global minimizer,
then the problem of finding such a minimizer is strongly stable
under perturbations in the model ($T$), in the data ($y$) and in the
regularization parameter ($\alpha$). This follows immediately from
the fact that \textit{(H2)} is stronger than \textit{(H2')}, the
relative compactness of BV-bounded sets in $L^p(\Omega)$ for
$p<\frac{d}{d-1}$ (see \cite{refb:Giusti-1984}) and Corollary
\ref{cor:strong-stability-Tikhonov}. For $p=\frac{d}{d-1}$ and $d\ge
2$ the problem is weakly stable, by virtue of Corollary
\ref{corollary:weak}.

\subsection{Penalization with powers of semi-norms associated to closed operators}\label{ssec:L}

\begin{thm}\label{teo:closed-operators-1}
Let $\X,\, \mathcal Z$ be reflexive Banach spaces, $\Y$ a normed
space, $T\in\cal{L}(\X,\Y)$ and $L:\D(L)\subset\X
\rightarrow\mathcal Z$ a closed linear operator such that the
range of $L$, $\mathcal R(L)$, is weakly closed. Assume further
that $T$ and $L$ are complemented, i.e. there exists a constant
$k>0$ such that $\|Tx\|^2+\norm{Lx}^2\geq k\norm{x}^2,\,\,
\forall\,\, x\in\D(L)$. Then, for any $q>1$, $\alpha>0$ and
$y\in\Y$ the functional
\begin{equation}
J_{\s L,\,q,\,\alpha}(x)\doteq \|Tx-y\|^2 + \alpha\|Lx\|^q, \qquad
x\in \D(L), \label{L-functional}
\end{equation}
has a unique global minimizer.
\end{thm}
\begin{proof}
Let $q>1$, $\D\doteq \D(L)$ and $W_{\s
L,q}:\D\longrightarrow\mathbb R^{\scriptscriptstyle
+}_{\scriptscriptstyle 0}$ defined by $W_{\s
L,q}(x)\doteq\norm{Lx}^q$. We will show that $T$ and $W_{\s L,q}$
satisfy the hypotheses \textit{(H1), (H2)} and \textit{(I3)}.
Hypothesis \textit{(H1)} is trivially satisfied since $W_{\s
L,q}(x)\ge 0\;\forall\; x\in \D$. To prove that \textit{(H2)}
holds, let $\{x_n\}\subset \D$ be a $W_{\s L,q}$-bounded sequence
such that $x_n\overset{w}{\rightarrow}x\in\D$. Then there exists a
constant $c<\infty$ such that $\norm{Lx_n}\leq
c\,\,\,\forall\,\,n\in\mathbb N$. Since the Banach space $\mathcal
Z$ is reflexive, there exist $z\in\mathcal Z$ and
$\{x_{n_j}\}\subset \{x_n\}$ such that
$Lx_{n_j}\overset{w}{\rightarrow} z$. Since $\R(L)$ is weakly
closed $z\in\R(L)$. Now, the operator $L^\dagger$, the
Moore-Penrose generalized inverse of $L$, is continuous (since
$\R(L)$ is closed), and therefore $P_{\s
\NU(L)^\perp}x_{n_j}=L^\dagger L x_{n_j}
\overset{w}{\longrightarrow}L^\dagger z$ (where
$P\s_{\NU(L)^\perp}$ is the orthogonal projection of $\X$ onto
$\NU(L)^\perp$). Since $x_{n_j}= P_{\s \NU(L)^\perp} x_{n_j} +
P_{\s \NU(L)}x_{n_j}$ it follows that $P_{\s
\NU(L)}x_{n_j}\overset{w}{\longrightarrow}x-L^\dagger z$ and
therefore $x-L^\dagger z \in \NU(L)$ (since $\NU(L)$ is weakly
closed, $L$ being closed). Hence $0=L(x-L^\dagger z)=Lx-LL^\dagger
z=Lx-P_{\s \overline{\R(L)}}z=Lx-z$. Thus $z=Lx$ and $W_{\s
L,q}(x)=\|Lx\|^q=\|z\|^q\le\liminf_{j\to\infty}\|Lx_{n_j}\|^q=\liminf_{j\to\infty}
W_{\s L,q}(x_{n_j})$, where the inequality follows from the fact
that $Lx_{n_j}\overset{w}{\longrightarrow}z$ and the weak lower
semicontinuity of the norm in $\Z$. This proves \textit{(H2)}.

To prove that \textit{(I3)} holds, let $\{x_n\}\subset \D$ be a
$T$-$W_{\s L,q}$-bounded sequence. By the complementation
condition it follows that $\{x_n\}$ is bounded in $\X$ and by the
reflexivity of $\X$ there must exist a subsequence
$\{x_{n_j}\}\subset \{x_n\}$ and $x\in \X$ such that
$x_{n_j}\overset{w}{\longrightarrow}x$. It only remains to be
proved that $x\in \D=\D(L)$. For that observe that since
$\{x_{n_j}\}$ is a $W_{\s L,q}$-bounded sequence such that
$x_{n_j}\overset{w}{\longrightarrow}x$, following the same steps
as in the proof of \textit{(H2)} above, we obtain that there
exists $z\in \R(L)$ such that $x-L^\dagger z\in \NU(L)$. Since
$L^\dagger z\in \NU(L)^\perp\subset \D(L)$ it then follows that
$x\in \D(L)$. This finally proves that \textit{(I3)} holds.

Now, since hypothesis \textit{(H2)} implies hypothesis
\textit{(I2)} (see Remark \ref{obs:H-weakers-1}), Theorem
\ref{teo:existence-weak} now implies that for any $\alpha>0$,
$y\in\Y$, the functional $J_{\s L,\,q,\,\alpha}(x)$ defined by
(\ref{L-functional}), has a global minimizer on $\D(L)$. Since
$q>1$, from the complementation condition it follows easily that
$J_{\s L,\, q,\,\alpha}$ is strictly convex and therefore such a
global minimizer is unique.
\end{proof}

\medskip
It is appropriate point out here that the above hypotheses on $L$
are satisfied by most differential operators and that the
complementation condition holds, for instance, whenever $dim\,
\NU(L)<\infty$ and $\NU(T)\cap\NU(L)=\{0\}$. Also, the previous
theorem provides existence for any $q>0$. However uniqueness can
only be guaranteed for $q>1$ and, if $T$ is injective, also for
$q=1$.
\medskip

The next Lemma shows that the problem of finding the global
minimum of (\ref{L-functional}) is weakly stable under
perturbations on $y$, $\alpha$ and $T$.
\begin{lem}\label{lemma:closed-operators-2}
Let $\X,\, \Y,\,\mathcal Z,\, T,\, L ,\, \D$ as in Theorem
\ref{teo:closed-operators-1}, $q>1$, $y,\;y_n\in\Y$,
$\alpha,\;\alpha_n\ge 0$, $T_n\in\cal{L}(\X,\Y)$, $n=1,2,\dots$,
 and $J_{\s L, q, \alpha}, J_n, n=1,2,...$, functionals
on $\D$ defined by
\begin{eqnarray}
J_{\s L,\,q,\,\alpha}(x)&\doteq\|Tx-y\|^2 + \alpha \norm{Lx}^q,\label{N1}\\
J_n(x)&\doteq\|T_nx-y_n\|^2 + \alpha_n \norm{Lx}^q\label{N2}.
\end{eqnarray}
Assume that $\alpha_n\to \alpha$, $y_n\to y$ as $n\to \infty$ and
that $T_nx\to Tx$ uniformly for $x$ in $L$-bounded sets (i.e.
$\{T_n\}$ is $L$-uniformly consistent for $T$). Let $\bar x$ be
the unique minimizer of $J_{\s L,\,q,\,\alpha}$ and $x_n$ a global
minimizer of $J_n$. Then $x_n\overset{w}{\rightarrow} \bar x$.
\end{lem}
\begin{proof}
Let $W_{\s L,q}:\D\longrightarrow\mathbb R^{\scriptscriptstyle
+}_{\scriptscriptstyle 0}$ defined by $W_{\s
L,q}(x)\doteq\norm{Lx}^q$. In Theorem \ref{teo:closed-operators-1}
we proved that $T$ and $W_{\s L,q}$ satisfy hypotheses
\textit{(H1)}, \textit{(I2)} and \textit{(I3)}. Since by
hypothesis $\alpha_n\to\alpha$, $y_n\to y$ and $\{T_n\}$ is $W_{\s
L,q}$-uniformly consistent for $T$, the Lemma follows immediately
from Corollary \ref{corollary:weak-2}.
\end{proof}

From the point of view of applications of the Tikhonov-Phillips
methods, the weak stability result established by the previous
Lemma, although important, could render insufficient. A strong
stability result, at least on the data $y$ is highly desired. In
the next Lemma we show that such a result can be obtain by
imposing an additional hypothesis to the operator $L$.

\begin{lem} \label{lemma:closed-operators-3}
Let $\X,\, \mathcal Z,\,\Y,\, T,\, T_n\, L ,\, \D,\, q,\, W_{\s
L,q}, \, y,\,y_n,\, \alpha,\,\alpha_n,\, J_{\s L, q, \alpha}$,
$\bar x$, $x_n$ and $\,J_n, n=1,2,...$ as in Lemma
\ref{lemma:closed-operators-2}. Assume further that
$T$-$L$-bounded sets are compact in $\X$. Then $x_n\to \bar x$.
\end{lem}
\begin{proof} In Theorem \ref{teo:closed-operators-1} we proved that
$T$ and $W_{\s L,q}$  satisfy hypotheses \textit{(H1)} and
\textit{(I2)}. Since hypothesis \textit{(I2)} implies hypothesis
\textit{(I2')} and the compactness of $T$-$L$-bounded sets implies
\textit{(I3)'}, the lemma then follows from Corollary
\ref{cor:strong-stability-Tikhonov-2}.
\end{proof}
\begin{rem}\label{remark:closed-operators-4}
If $q=2$, under the same hypotheses of Lemma 4.2 one can get
continuity of the solutions of (\ref{N1}) with respect to $\alpha$
and $y$. This can be easily verified from the fact that the unique
global minimizer of (\ref{N1}) is given by $\bar x=(\alpha L^\ast
L+T^\ast T)^{-1}T^\ast y$. Thus, if $x_n$ is the minimizer of
(\ref{N2}) with $T_n=T \;\,\forall\, n$, then one has that
\begin{equation}\label{N3}
\bar x-x_n = (\alpha-\alpha_n)\left(\alpha L^\ast L + T^\ast
T\right)^{-1} L^\ast L\, x_n  + \left(\alpha L^\ast L + T^\ast
T\right)^{-1}T^\ast(y-y_n).
\end{equation}
Suppose now that $\alpha_n\to\alpha$ and $y_n\to y$. Then by Lemma
\ref{lemma:closed-operators-2} $x_n\overset{w}{\rightarrow} \bar
x$ and therefore $\{x_n\}$ is bounded. Also, since
$\|Tx\|^2+\norm{Lx}^2\geq k\norm{x}^2$ it follows that the
operators $\left(\alpha L^\ast L + T^\ast T\right)^{-1} L^\ast L$
and $\left(\alpha L^\ast L + T^\ast T\right)^{-1}T^\ast$ are both
bounded. In fact $\left(\alpha L^\ast L + T^\ast T\right)^{-1}
L^\ast L\le \alpha^{-1}I$ and $\left(\alpha L^\ast L + T^\ast
T\right)^{-1}$ $\le \frac{1}{k\min(\alpha,1)}$. Hence, it follows
from (\ref{N3}) that $x_n\to\bar x$.
\end{rem}

\bigskip
%
%
\subsection{Penalization by linear combination of powers of semi-norms associated to closed operators}
\label{ssec:comb lineal} We study here the case of generalized
Tikhonov-Phillips regularization methods for which the functional
$W(\cdot)$ in (\ref{eq:prob opt}) is of the form
$W(x)\doteq\sum_{i=1}^N \alpha_i \|L_ix\|^{q_i}$, where the
$L_i$'s are closed operators. We start with the main existence and
uniqueness result.
\begin{thm}\label{teo:hybrid-1}
Let $\X,\,\mathcal Z_1,\,\mathcal Z_2,\ldots, \mathcal Z_N$ be
reflexive Banach spaces, $\Y$ a normed space,
$T\in\cal{L}(\X,\Y)$, $\mathcal D$ a subspace of $\X$, $L_i:
\mathcal D \longrightarrow \mathcal Z_i,$ $ i=1,2,...,N,$ closed
linear operators with $\R(L_i)$ weakly closed for every $1\leq
i\leq N$ and such that $T,L_1,L_2,\dots,L_N$ are complemented,
i.e. there exists a constant $k>0$ such that $\|Tx\|^2
+\sum_{i=1}^N\|L_ix\|^2\ge k\|x\|^2$, $\forall\,x\in \D$. Then,
for any $y\in \Y$,
$\alpha_1,\alpha_2,\dots,\alpha_N\in\mathbb{R}^+$ and
$q_1,q_2,\dots,q_N\in\mathbb{R}$, $q_i>1\;\forall\,i=1,2,\dots,N$,
the functional
\begin{equation}\label{eq:N4}
J(x)\doteq \norm{Tx-y}^2+\sum_{i=1}^N\alpha_i\norm{L_ix}^{q_i},
\end{equation}
has a unique global minimizer.
\end{thm}
\begin{proof}
Let $y\in \Y$, $\alpha_i>0$, $q_i>1$, $i=1,2,\dots,N$ and define
$\vec\alpha\doteq(\alpha_1,\alpha_2,\dots,\alpha_N)^T,\, \vec
q\doteq(q_1,q_2,\dots,q_N)^T$, the normed space $\Z\doteq
\bigotimes_{i=1}^N \Z_i$, $\vec L:\X\to\Z$ as $\vec L
x\doteq(L_1x,L_2x,\dots,L_Nx)^T$, and the functional $W_{\vec
L,\vec q,\vec\alpha}:\D\to\mathbb{R}^+_0$ by $W_{\vec L,\vec
q,\vec\alpha}(x)=\sum_{i=1}^N\alpha_i\norm{L_ix}^{q_i}$, so that
$J(x)=\|Tx-y\|^2+W_{\vec L,\vec q,\vec\alpha}(x)$. We will prove
that $T$ and $W_{\vec L,\vec q,\vec\alpha}$  satisfy the
hypotheses \textit{(H1)}, \textit{(H2)} and \textit{(I3)}. In
fact,  \textit{(H1)} is trivial and for \textit{(H2)}, let
$\{x_n\}\subset\D$ be a $W_{\vec L,\vec q,\vec\alpha}$-bounded
sequence such that $x_n\overset{w}{\to}x\in\D$. Then for every
$i=1,2,\dots,N$, the sequence $\{L_ix_n\}_{n=1}^\infty$ is bounded
in $\Z_i$ and since $\Z_i$ is reflexive there exist a subsequence
$\{x_{n_k}\}$ and $z_i\in\Z_i$ such that
$L_ix_{n_k}\overset{w}{\to}z_i$ as $k\to\infty$. Since $\R(L_i)$
is weakly closed, $z_i\in \R(L_i)$. By taken subsequences, we may
assume that such a subsequence is the same for all $i$, i.e.
$L_ix_{n_k}\overset{w}{\to}z_i$ as $k\to\infty$ for every
$i=1,2,\dots,N$.

Now, since $\R(L_i)$ is closed, $L_i^\dagger$ is bounded and
therefore $L_i^\dagger L_ix_{n_k}\overset{w}{\to}L_i^\dagger z_i$,
as $k\to\infty$, for all $i=1,2,\dots,N$. Since $L_i^\dagger L_i
=P_{\mathcal{N}(L_i)^\perp}$ is the orthogonal projection of $\X$
onto $\mathcal{N}(L_i)^\perp$, writing
$x_{n_k}=P_{\mathcal{N}(L_i)^\perp}x_{n_k} +
P_{\mathcal{N}(L_i)}x_{n_k}$, it follows that
$P_{\mathcal{N}(L_i)}x_{n_k}\overset{w}{\to} x-L_i^\dagger z_i$ as
$k\to\infty$ and therefore $x-L_i^\dagger z_i\in \mathcal{N}(L_i)$
(being $\mathcal{N}(L_i)$ closed, since $L_i$ is closed). Hence
for all $i=1,2,\dots,N$, it follows that $0=L_i(x-L_i^\dagger
z_i)=L_ix- P_{\overline{\R(L_i)}}z_i=L_ix-z_i$ (where the last
equality follows since $z_i\in \R(L_i)$\,). Thus,
$z_i=L_ix\;\forall\, i=1,2,\dots,N$. Then
\begin{equation*}
\|L_ix\|^{q_i}=\|z_i\|^{q_i}\le\liminf_{k\to\infty}
\|L_ix_{n_k}\|^{q_i},
\end{equation*}
(where the last inequality follows from the fact that $L_i x_{n_k}
\overset{w}{\to}z_i$ as $k\to\infty$ and the weak lower
semicontinuity of the norm in $\Z_i$), and therefore
\begin{eqnarray*}
W_{\vec L,\vec
q,\vec\alpha}(x) &=&\sum_{i=1}^N\alpha_i\|L_ix\|^{q_i}
 \le\sum_{i=1}^N\alpha_i\liminf_{k\to\infty}\|L_ix_{n_k}\|^{q_i}\\
 &\le& \liminf_{k\to\infty} \sum_{i=1}^N\alpha_i\|L_ix_{n_k}\|^{q_i}
 = \liminf_{k\to\infty}W_{\vec L,\vec q,\vec\alpha}(x_{n_k}).
\end{eqnarray*}
Thus \textit{(H2)} holds. That \textit{(I3)} also holds follows
from the complementation condition and the reflexivity of $\X$,
following the same steps as in Theorem
\ref{teo:closed-operators-1}. Since \textit{(H2)} implies
\textit{(I2)}, it now follows from Theorem
\ref{teo:existence-weak} that the functional $J(x)$ in
(\ref{eq:N4}) has a global minimizer on $\D$. Moreover, since
$q_i>1$ for all $i$, it follows from the complementation condition
that $J(\cdot)$ is strictly convex and therefore such a minimizer
is unique.
\end{proof}
Under the same hypotheses of Theorem \ref{teo:hybrid-1} one has
that the solution of (\ref{eq:N4}) is weakly stable under
perturbations in the data $y$, in the parameters $\alpha_i$ and in
the model operator $T$. More precisely we have the following
result.
\begin{lem}\label{lemma:hybrid-2}
Let all the hypotheses of Theorem \ref{teo:hybrid-1} hold. Let
also $y,y_n\in\Y$, $T_n\in\cal{L}(\X,\Y)$, $n=1,2,\dots$, such
that $y_n\to y$, $\{T_n\}$ is $\vec L$-uniformly consistent for
$T$ and for each $i=1,2,\dots,N$, let
$\{\alpha_i^n\}_{n=1}^\infty\subset\mathbb{R}^+$ such that
$\alpha_i^n\to\alpha_i$ as $n\to\infty$. If $x_n$ is a global
minimizer of the functional
\begin{equation}\label{eq:N5}
J_n(x)\doteq\|T_nx-y_n\|^2+\sum_{i=1}^N\alpha_i^n\|L_ix\|^{q_i},
\end{equation}
then $x_n\overset{w}{\to}\bar x$, where $\bar x$ is the unique
minimizer of (\ref{eq:N4}).
\end{lem}
\begin{proof}
Let $W\doteq W_{\vec L,\vec q,\vec\alpha}$ as in Theorem
\ref{teo:hybrid-1}. From the hypotheses it follows easily that
$\{J_n\}$ is $T$-$W$-coercive and $W$-uniformly consistent for
$J$.

Let $x_n$ be the unique minimizer of $J_n$. Then $J_n(x_n)\le
J_n(\bar x)$, $\forall\, n$. Therefore
\begin{equation}\label{eq:P0}
\limsup_{n\to\infty}J_n(x_n)\le \limsup_{n\to\infty}J_n(\bar
x)=J(\bar x)<\infty,
\end{equation}
where the equality follows from the $W$-uniform consistency of
$\{J_n\}$ for $J$. But since $\{J_n\}$ is $T$-$W$-coercive it then
follows that $\{x_n\}$ is $T$-$W$-bounded. We claim that
$x_n\overset{w}{\to}\bar x$. In fact, suppose that is not the
case. Then, there exists a subsequence $\{x_{n_j}\}\subset\{x_n\}$
such that no subsequence of $\{x_{n_j}\}$ converges weakly to
$\bar x$. But since  $\{x_{n_j}\}$ is $T$-$W$-bounded and $\X$ is
reflexive, there exist $x^\ast\ne \bar x$ and
$\{x_{n_{j_k}}\}\subset\{x_{n_j}\}$ such that
$x_{n_{j_k}}\overset{w}{\to}x^\ast$. Following the same steps as
in Theorem \ref{teo:hybrid-1} we obtain that there exists a
subsequence $\{x_{n_{j_{k_\ell}}}\}\subset\{x_{n_{j_k}}\}$ and
$z_i\in \Z_i$, $i=1,2,\dots,N$, such that $L_ix_{n_{j_{k_\ell}}}
\overset{w}{\to}z_i=L_ix^\ast$ as $\ell\to\infty$,
$\forall\,i=1,2,\dots,N$, and
\begin{equation}\label{eq:P1}
W(x^\ast)\le\liminf_{\ell\to\infty} W\left(x_{n_{j_{k_\ell}}}
\right).
\end{equation}
Also, since $\{x_{n_{j_{k_\ell}}}\}$ is $W$-bounded and $\{J_n\}$
is $W$-uniformly consistent for $J$, it follows that
\begin{equation}\label{eq:P2}
\lim_{\ell\to\infty}\left(J\left(x_{n_{j_{k_\ell}}}\right)-
J_{n_{j_{k_\ell}}}\left(x_{n_{j_{k_\ell}}}\right) \right)\;=\;0.
\end{equation}
Hence 
\begin{eqnarray*}
J(x^\ast)&=\|Tx^\ast-y\|^2+W(x^\ast)\\
&\le\liminf_{\ell\to\infty}\|Tx_{n_{j_{k_\ell}}}  -y\|^2 +
\liminf_{\ell\to\infty}W(x_{n_{j_{k_\ell}}})\qquad\qquad\qquad\text{(by
(\ref{eq:P1})\;)}\\
&\le\liminf_{\ell\to\infty}\left(\|Tx_{n_{j_{k_\ell}}} -y\|^2 +
W(x_{n_{j_{k_\ell}}}) \right)\\
&= \liminf_{\ell\rightarrow\infty} J\left(x_{n_{j_{k_\ell}}}\right) \\
&\leq \limsup_{\ell\rightarrow\infty} J\left(x_{n_{j_{k_\ell}}}\right) \\
&=  \limsup_{\ell\rightarrow\infty}\left[\left(
J\left(x_{n_{j_{k_\ell}}}\right)-J_{n_{j_{k_\ell}}}
\left(x_{n_{j_{k_\ell}}}\right) \right) +J_{n_{j_{k_\ell}}}\left(x_{n_{j_{k_\ell}}}\right)\right]\\
&\leq
\limsup_{\ell\rightarrow\infty}\left(J\left(x_{n_{j_{k_\ell}}}\right)-J_{n_{j_{k_\ell}}}\left(x_{n_{j_{k_\ell}}}\right)
\right)+\limsup_{\ell\rightarrow\infty} J_{n_{j_{k_\ell}}}\left(x_{n_{j_{k_\ell}}}\right)\\
&=\limsup_{\ell\rightarrow\infty} J_{n_{j_{k_\ell}}}\left(x_{n_{j_{k_\ell}}}\right) \qquad\qquad\text{(by (\ref{eq:P2})\,)}\\
&\leq J(\bar x) \quad\qquad\qquad\qquad\qquad\qquad\text{(by
(\ref{eq:P0}) since $\{x_{n_{j_{k_\ell}}}\}\subset\{x_n\}$\,)}
\end{eqnarray*}
Since $\bar x$ is the unique minimizer of $J$ it would then follow
that $x^\ast=\bar x$, contradicting our previous result that
$x^\ast\ne\bar x$. This contradiction came from the assumption
that $x_n$ did not converge weakly to $\bar x$. Hence
$x_n\overset{w}{\to}\bar x$.
\end{proof}
\begin{lem}\label{lemma:hybrid-3} Under the same hypotheses of
Lemma \ref{lemma:hybrid-2}, if\; $T$-$\vec L$-bounded sets are
compact in $\D$, then strong stability holds, i.e., $x_n\to \bar
x$.
\end{lem}
\begin{proof}
Let $x_n$ denote the global minimizer of $J_n$ and $W=W_{\vec L
,\vec q,\vec\alpha}$. In Lemma \ref{lemma:hybrid-2} it was proved
that the sequence $\{x_n\}$ is $T$-$W$-bounded. Suppose that $x_n
\nrightarrow \bar x$. Then there exists a subsequence
$\{x_{n_j}\}\subset\{x_n\}$ such that no subsequence of
$\{x_{n_j}\}$ converges to $\bar x$. But since $\{x_{n_j}\}$ is
$T$-$W$-bounded, now by compactness hypothesis there must exist
$x^\ast\in \D$, $x^\ast\ne \bar x$, and a subsequence
$\{x_{n_{j_k}}\}\subset\{x_{n_j}\}$ such that $x_{n_{j_k}}\to
x^\ast$ as $k\to\infty$. Using the $W$-uniform consistency of
$\{J_n\}$ for $J$ and following similar steps as in Lemma
\ref{lemma:hybrid-2} one obtains that $J(x^\ast)\le J(\bar x)$.
Since $\bar x$ is the unique minimizer of $J$ it would then follow
that $x^\ast=\bar x$, contradicting our previous result that
$x^\ast\ne\bar x$. Therefore we must have that $x_n\to\bar x$.
\end{proof}
\begin{rem} Here again, for the case $q_i=2\;\forall\, i$, strong
continuity of the solution of the functional $J(x)$ in
(\ref{eq:N4}) with respect to the data $y$ and the parameters
$\alpha_i$ follow without any further hypotheses than those in
Lemma \ref{lemma:hybrid-2}. This result follows easily from the
fact that in such a case the unique global minimizer of
(\ref{eq:N4}) is given by $\bar x=\left( T^\ast
T+\sum_{i=1}^N\alpha_i L_i^\ast L_i\right)^{-1}T^\ast y$. Thus, if
$x_n$ is the minimizer of (\ref{eq:N5}) with $T_n=T\;\forall\, n$,
then one has that
\begin{eqnarray}\label{eq:Q1}
\bar x-x_n &=& \left( T^\ast T+\sum_{i=1}^N\alpha_i L_i^\ast
L_i\right)^{-1} \sum_{i=1}^N(\alpha_i^n-\alpha_i)L_i^\ast L_i
\,x_n\nonumber\\
&&\;\;+ \left( T^\ast T+\sum_{i=1}^N\alpha_i L_i^\ast
L_i\right)^{-1} T^\ast (y-y_n).
\end{eqnarray}
Now, from the complementation condition
$\|Tx\|^2+\sum_{i=1}^N\|L_ix\|^2\ge k\|x\|^2,\forall\, x\in\D$, it
follows easily that
\begin{equation}\label{eq:Q2}
0\le \left(T^\ast T+\sum_{i=1}^N\alpha_i L_i^\ast
L_i\right)^{-1}\le
 \frac{1}{k\min\left(1,\min_{1\le i\le N}\alpha_i\right)},
\end{equation}
and also
\begin{eqnarray}\label{eq:Q3}
&\left\|\left(T^\ast T+\sum_{i=1}^N\alpha_i L_i^\ast
L_i\right)^{-1} \sum_{i=1}^N(\alpha_i^n-\alpha_i)L_i^\ast L_i\,x
\right\| \nonumber\\
&\quad\le\frac{\max_{1\le i\le N}|\alpha_i^n-\alpha_i|}{\min_{1\le
i\le N}\alpha_i}\,\|x\|,\quad\forall\, x\in \D.
\end{eqnarray}
Using (\ref{eq:Q3}) and (\ref{eq:Q2}) in (\ref{eq:Q1}) we obtain
that
\begin{equation}\label{eq:Q4}
\|\bar x-x_n\|\le \frac{\max_{1\le i\le
N}|\alpha_i^n-\alpha_i|}{\min_{1\le i\le N}\alpha_i}\,\|x_n\| +
\frac{\|T^\ast\|}{k\min\left(1,\min_{1\le i\le
N}\alpha_i\right)}\,\|y-y_n\|.
\end{equation}
Now since by Lemma \ref{lemma:hybrid-2} $x_n\overset{w}{\to}\bar
x$, it follows that $\{x_n\}$ is bounded. Since $y_n\to y$ and
$\alpha_i^n\to\alpha_i\;\forall\, i=1,2,\dots,N$, as $n\to\infty$,
it finally follows from (\ref{eq:Q4}) that $x_n\to x$.
\end{rem}


\section{Applications to Image Restoration} \label{sec:numerical-examples}

The purpose of this section is to present an application to a
simple image restoration problem. The main objective is to show
how the choice of the penalizer in a generalized Tikhonov-Phillips
functional can affect the reconstructed image.

The basic mathematical model for image blurring is given by the following Fredholm integral
equation
\begin{equation}\label{eq:restau}
K\,f(x,y)\doteq\int\int_{\Omega}
k(x,y,x',y')f(x',y')dx'dy'=g(x,y),
\end{equation}
where $\Omega\subset\mathbb R^{2}$ is a bounded domain, $f\in\X\doteq L^{2}(\Omega)$
represents the original image and $k$ is the so called ``point spread function'' (PSF).
For the examples shown below we used a PSF of ``atmospheric turbulence'' type
\begin{eqnarray}\label{eq:PSF en aplicaciones}
k(x,y,x',y')&=\frac{\kappa}{\pi}\exp\left(-\kappa\norm{(x,y)-(x',y')}^2\right),
\end{eqnarray}
with $\kappa=6$. It is well known that with this PSF the operator
$K$ in (\ref{eq:restau}) is compact with infinite dimensional
range and therefore $K^{\dag}$, the Moore-Penrose inverse of $K$,
is unbounded.

Generalized Tikhonov-Phillips methods with different penalizers where used to obtain
regularized solutions of the problem
\begin{equation}\label{eq:PI}
K\,f=g.
\end{equation}
The data $g$ was contaminated with a $1\%$ zero mean Gaussian
noise (i.e. standard deviation of the order of $1\%$ of
$\norm{g}_{\s\infty}$). Minimizers of functionals of the form
\begin{equation} \label{eq:general-functional}
J_\alpha(f)=\|Kf-\tilde g\|^2 +\alpha\, W(f)
\end{equation}
were found for different penalizers $ W(f)$, where $\tilde g$
represents the noisy version of $g$. In all cases the value of the
regularization parameter $\alpha$ was approximated by using the
L-curve method (\cite{refb:Engl-Hanke-Neubauer-1996},
\cite{ref:Hansen-1992}, \cite{ref:Hansen-O'Leary-1993}).

\begin{figure}[H]
\begin{center}
\subfigure[Original image (unknown). \label{fig:original}
 ]{\includegraphics[width=2in]{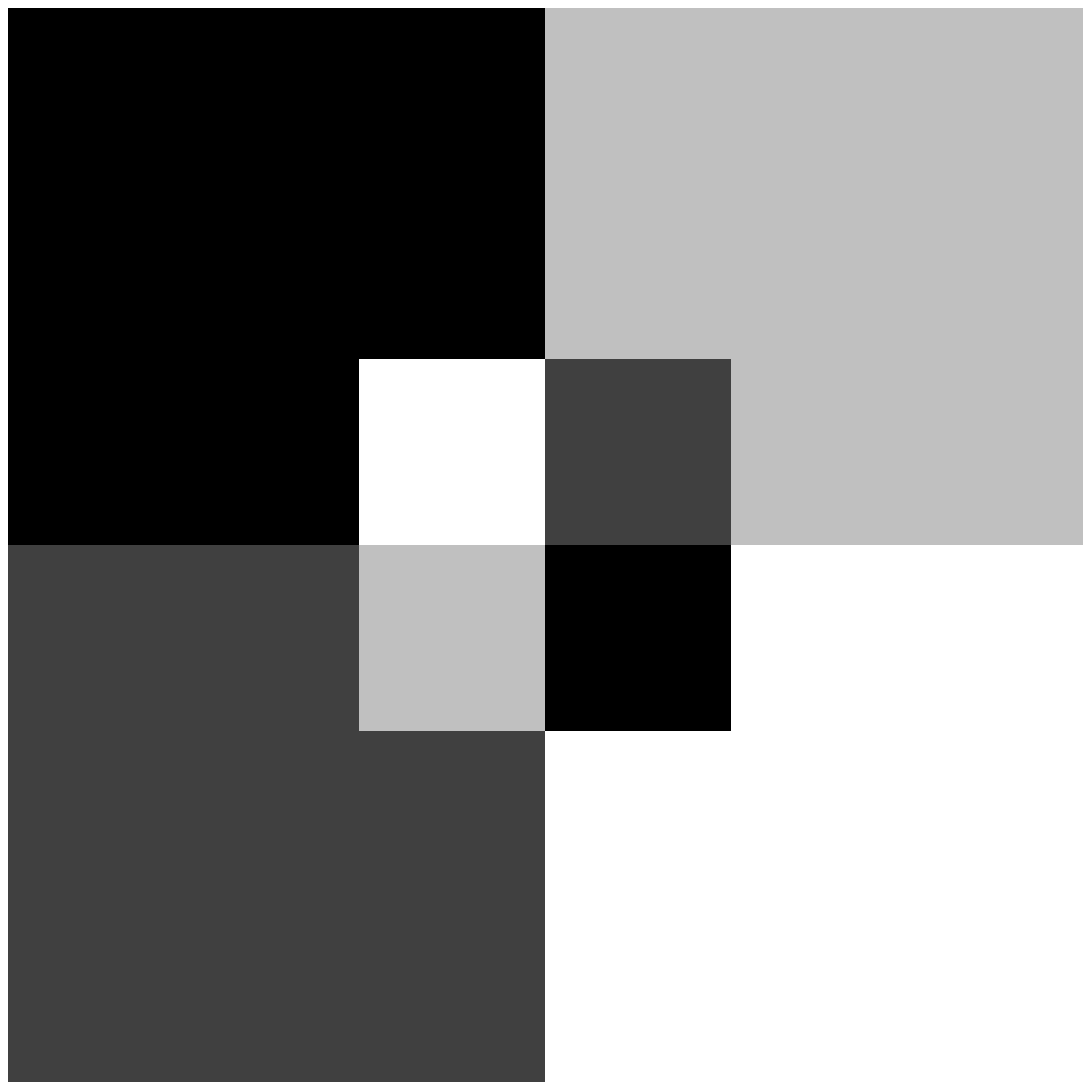}}\quad
 \subfigure[Blurred noisy image (data).
 \label{fig:blurred}]{\includegraphics[width=2in]{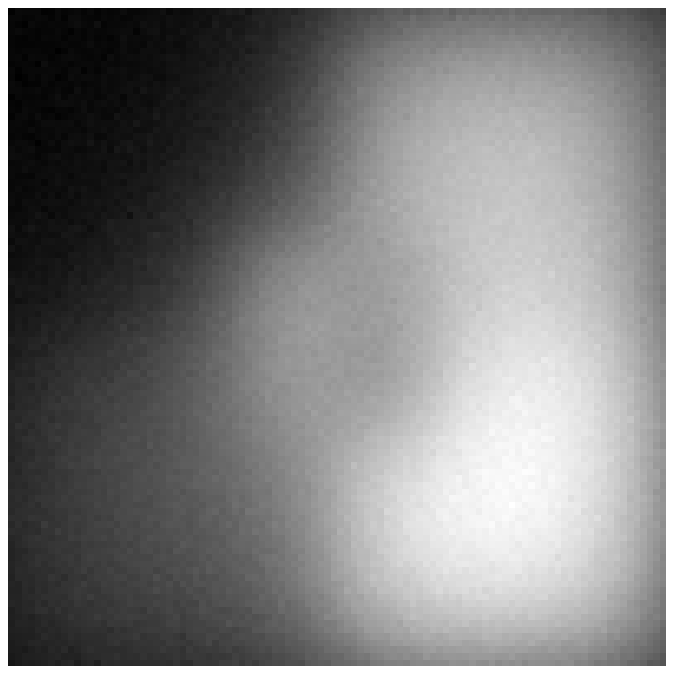}}

 \subfigure[Tikhonov-Phillips of order zero, $W(f)=\|f\|^2$.
  \label{fig:tik0}]{\includegraphics[width=2in]{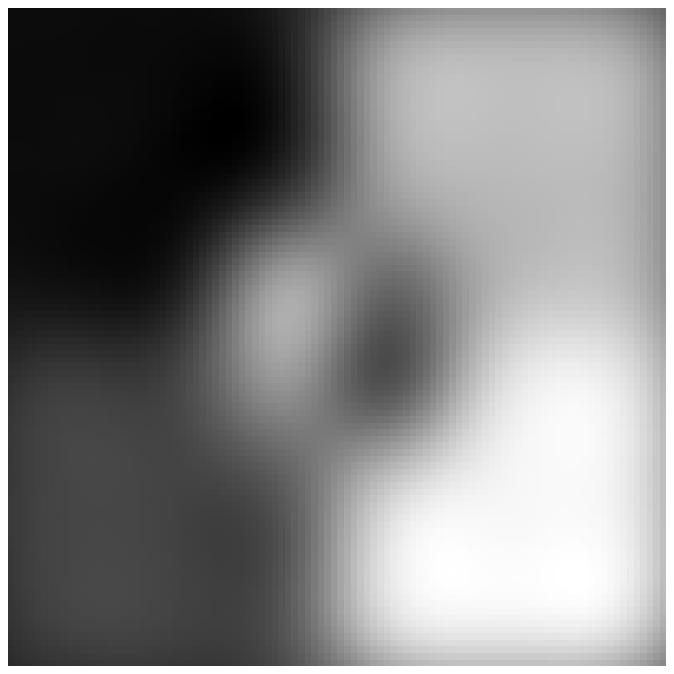}}\quad
 \subfigure[Tikhonov-Phillips of order one, $W(f)=\norm{\nabla
 f}^2$.
 \label{fig:tik1}]{\includegraphics[width=2in]{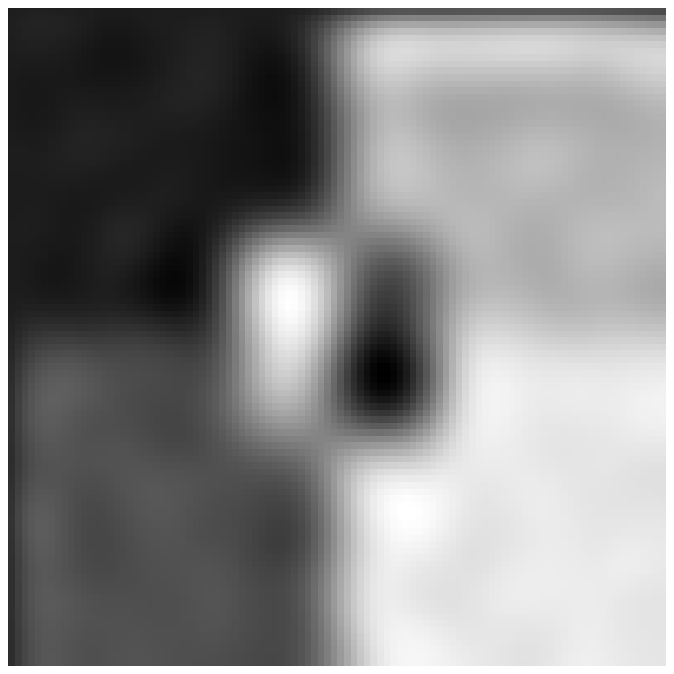}}
 \end{center}
\caption{Original image (a), blurred noisy image (b) and
regularized solutions obtained with the classical
Tikhonov-Phillips methods of order zero (c) and one (d).
}\label{fig:one}
\end{figure}

Figures \ref{fig:original} and \ref{fig:blurred} show the original
image (unknown in real life problems) and the blurred noisy image
which constitutes the data for the inverse problems, respectively.
Figures \ref{fig:tik0} and \ref{fig:tik1} show the reconstructions
obtained with the classical Tikhonov-Phillips methods of order
zero and one, corresponding to $W(f)=\|f\|^2$ and $W(f)=\|\nabla
f\|^2$, respectively.

\begin{figure}[H]
\begin{center}
\subfigure[The curve $\gamma$ providing the structural
information. \label{fig:infostruct}
]{\includegraphics[width=2in]{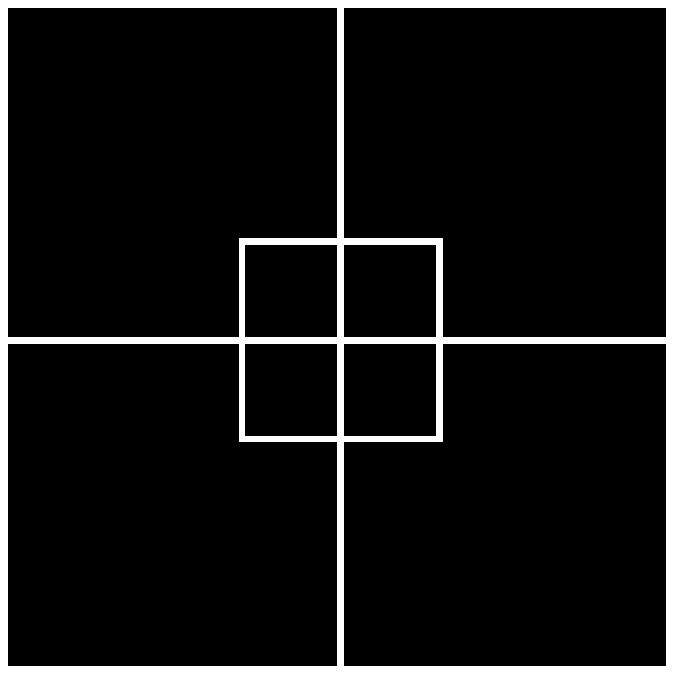}}\quad
\subfigure[Structural penalizer $W(f)=\|Lf\|$,  $c=5$.
\label{fig:structc5}]{\includegraphics[width=2in]{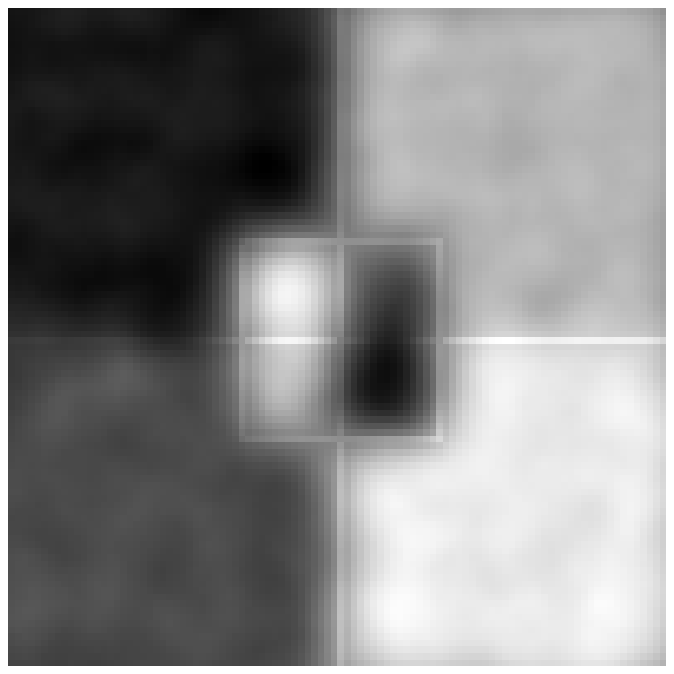}}

\subfigure[Structural penalizer $W(f)=\|Lf\|$,  $c=20$. \label{fig:structc20}
 ]{\includegraphics[width=2in]{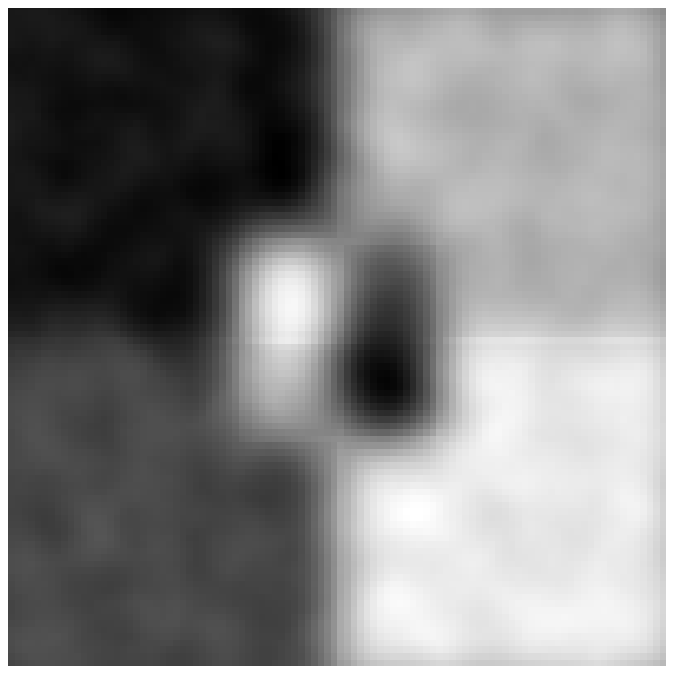}}
\quad
\end{center}
\caption{Structural information (a), reconstructed image with
structural information penalizer and $c=5$ (b) and $c=20$
(c)}\label{fig:two}
\end{figure}

Figures \ref{fig:structc5} and \ref{fig:structc20} show the
reconstructions obtained with a structural information penalizer
of the form $W(f)=\|Lf\|^2$ where the operator $L$ is constructed
as in \cite{ref:Kaipio-Kolehmainen-Vauhkonen-1999}, including the
information of the curve $\gamma$ depicted in Figure
\ref{fig:infostruct}, where it is expected that the original image
have steep gradients. The operator $L$ is constructed so as to
capture this structural prior information. The discretization of
$L$ is given by $\int_\Omega\|A(x)\nabla f(x)\|^2\;dx$ with
$A(x)=I-\left( 1+c\|\nabla
\gamma(x)\|^2\right)^{-1}\nabla\gamma(x)\left(\nabla\gamma(x)\right)^T$,
where $c$ is a positive constant. In this way, if $\|\nabla
\gamma(x)\|$ is large, the functional $W(f)$ penalizes only very
mildly all intensity changes occurring in the direction of
$\nabla\gamma(x)$ (see
\cite{ref:Kaipio-Kolehmainen-Vauhkonen-1999} for more details).

\begin{figure}[H]
\begin{center}
 \subfigure[\label{fig:hibrido1}
 ]{\includegraphics[width=2in]{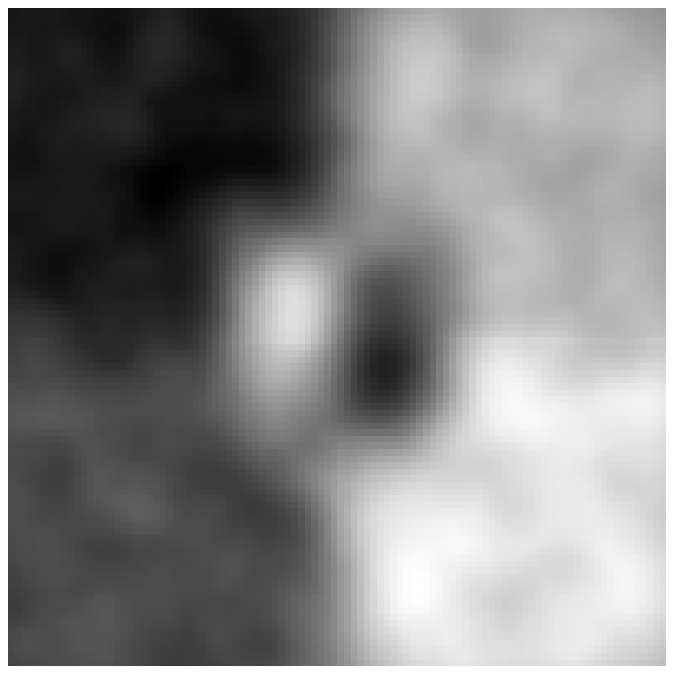}}
 \subfigure[\label{fig:hibrido2}
 ]{\includegraphics[width=2in]{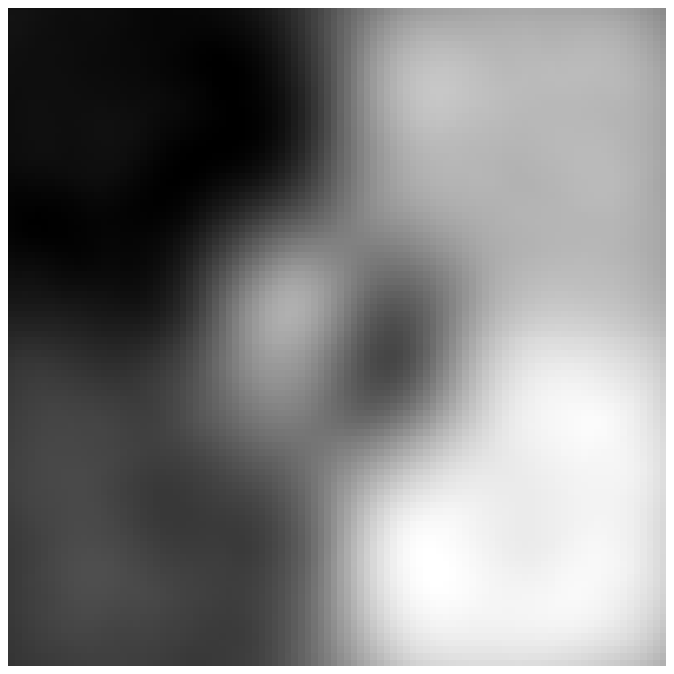}}
  \end{center}
\caption{Recostructed images with hybrid penalizers:
$W(f)=\frac45\|f\|^2+\frac15\|Lf\|^2$; $c=5$ (a) and $c=20$
(b)}\label{fig:three}
\end{figure}

Figures \ref{fig:hibrido1} and \ref{fig:hibrido2} correspond to
images reconstructed with hybrid Tikhonov-structural penalizers
$W(f)=\frac45\|f\|^2+\frac15\|Lf\|^2$, with $c=5$ and $c=20$,
respectively.

A comparison of the images obtained with the different methods
clearly show that the choice of the penalizer in Tikhonov-Phillips
method can greatly affect the obtained approximated solution. In
this particular case we observe how the classical order-zero
method tends to smooth out boundaries and edges and, while the
order-one method does a better job, the inclusion of the
structural information through the operator $L$ results in a
significant improvement. Although the main objective of this
article is theoretical in nature, providing sufficient conditions
on the model operators and the penalizers for the existence,
uniqueness and stability of solutions of the corresponding
generalized Tikhonov-Phillips methods, the previous applications
to image restoration were included to better emphasize the
importance of the adequate choice of the penalizer.

\section{Conclusions}
In this article sufficient conditions on the penalizers in
generalized Tikhonov-Phillips functionals guaranteeing existence,
uniqueness and stability of the minimizers where found. The
particular cases in which the penalizers are given by the bounded
variation norm, by powers of seminorms and by linear combinations
of powers of seminorms associated to closed operators, were
studied. Several examples were presented and a few results on
image restoration were shown to illustrate how the choice of the
penalizer can greatly affect the regularized solutions.

\ack This work was supported in part by Consejo Nacional de
Investigaciones Cient\'{\i}ficas y T\'{e}cnicas, CONICET, through PIP
2010-2012 Nro. 0219, by Universidad Nacional del Litoral, U.N.L.,
through project CAI+D 2009-PI-62-315, by Agencia Nacional de
Promoci\'{o}n Cient\'{\i}fica y tecnol\'{o}gia ANPCyT, through project PICT
2008-1301 and by the Air Force Office of Scientific Research,
AFOSR, through Grant FA9550-10-1-0018.

\section*{References}
\bibliographystyle{amsplain}
\bibliography{ref}

\end{document}